\documentclass[a4paper]{article}
\usepackage{amssymb}
\usepackage{amsthm}
\usepackage{amsmath}
\usepackage{latexsym}
\usepackage[T1]{fontenc}
\usepackage[latin1]{inputenc}

\def \R{\mathbb{R}}

\def \N{\mathbb{N}}
\def \E{\mathbb{E}}

\def \L{\mathbb{L}}

\def \bf{\textbf}
\def \it{\textit}

\def \bop {\noindent\textbf{Proof }}
\def \eop {\hbox{}\nobreak\hfill
\vrule width 2mm height 2mm depth 0mm
\par \goodbreak \smallskip}
\def \ni {\noindent}

\def \F{\mathcal{F}}

\def \eop {\hbox{}\nobreak\hfill \vrule width 2.0mm height 1.8mm depth 0mm
\par \goodbreak \smallskip}

\numberwithin{equation}{section}
\newcommand{\dint}{\displaystyle\int}
\newcommand{\1}{1\!\!1}
\newtheorem{definition}{Definition}[section]
\newtheorem{theorem}{Theorem}[section]
\newtheorem{proposition}{Proposition}[section]
\newtheorem{lemma}{Lemma}[section]
\newtheorem{remark}{Remark}[section]

\def \R{\mathbb{R}}

\def \N{\mathbb{N}}
\def \E{\mathbb{E}}

\def \L{\mathbb{L}}

\def \bf{\textbf}
\def \it{\textit}

\def \bop {\noindent\textbf{Proof }}
\def \eop {\hbox{}\nobreak\hfill
\vrule width 2mm height 2mm depth 0mm
\par \goodbreak \smallskip}
\def \ni {\noindent}

\def \F{\mathcal{F}}

\begin{document}
\title{$p$-integrable solutions to multidimensional BSDEs
and degenerate systems of PDEs with logarithmic
nonlinearities.}
 \author{\quad K. Bahlali$^1$\; E. H. Essaky$^2$\; M. Hassani$^2$}
\date{}
 \maketitle
\begin{center}
{$^{1}$ IMATH, UFR Sciences, USVT, B.P. 132, 83957 La Garde Cedex,
France. \\ \small{e-mail: bahlali@univ-tln.fr}}
\\
{$^{2}$ Universit\'{e} Cadi Ayyad Facult\'{e} Poly-disciplinaire
D\'{e}partement de
\\ Math\'{e}matiques et d'Informatique, B.P. 4162 Safi, Maroc. \\ \small{e-mail:
essaky@ucam.ac.ma, medhassani@ucam.ac.ma}}
\end{center}

\begin{abstract}
 We study multidimensional backward stochastic
 differential equations
(BSDEs) which cover the logarithmic nonlinearity $u\log u$.
 More precisely, we establish the existence and uniqueness as well as the stability
 of $p$-integrable solutions $(p>1)$ to multidimensional
BSDEs with a $p$-integrable terminal condition and a super-linear
growth  generator in the  both variables $y$ and $z$. This is done
with
 a generator $f(y,z)$ which can be neither locally monotone in the
variable $y$ nor locally Lipschitz in the variable $z$. Moreover,
 it is not uniformly continuous.
As application, we establish the existence and uniqueness of Sobolev
solutions to possibly degenerate systems of semilinear parabolic
PDEs with super-linear growth generator and an $p$-integrable
terminal data. Our result cover, for instance, certain (systems of)
 PDEs
 arising in physics.
\end{abstract}


\section{\bf{Introduction}}

The logarithmic nonlinearity $u\log u$ appears in certain
differential equations arising in physics (see e. g. \cite{BM1,BM2,
CH, H,  PZ}) and in the theory of continuous-state branching
processes (see e. g. \cite{BeLeg, FS, FW}). For instance, the Cauchy problem
\begin{equation}\label{PDEFleishman}
\left\{
\begin{aligned}
&\dfrac{\partial u}{\partial t} - \Delta u + u\log u =  0   \ \ \hbox{on} \ \ (0,\ \infty)\times \R^d  \\
&u(0^+) = \varphi >0
\end{aligned}
\right.
\end{equation} is related to super processes with Neveu's branching mechanism, see e. g \cite{FS}. On the other hand, the logarithmic nonlinearity is
also interesting in itself since it is neither  locally monotone
 nor uniformly continuous. In this paper, we give a BSDEs approach which allows to cover this kind of nonlinearity.

Let $(W_t)_{0\;\;\leq\;\; t \leq T}$ be a $r$-dimensional Wiener process
defined on a complete probability space $(\Omega,\F, P)$.
$(\F_t)_{0\;\;\leq\;\; t \leq T}$ denote the natural filtration of $(W_t)$ such
that $\F_0$ contains all P-null sets of $\F$, and $\xi$ be an $\F_T$-measurable
$d$-dimensional random variable. Let $f$ be an $\R^d$-valued function defined
on $[0,T]\times\Omega\times\R^d\times\R^{d\times r}$ such that for every
$(y,z)\in \R^d\times\R^{d\times r}$, the map $(t,\omega)\longrightarrow
f(t,\omega,y,z)$ is $\F_t$-progressively measurable. The BSDEs under consideration
are of the form,
 $$  Y_t = \xi + \dint_t^T f(s,Y_s,Z_s)ds - \dint_t^T Z_s dW_s \quad\quad
0\leq t\leq T \leqno (E^{(\xi,f)}) $$
The data $\xi$ and $f$ are respectively called  the terminal condition and the
coefficient or generator.
\par
The present paper is a developpment of \cite{BEH}, and it constitute
a natural continuation of our previous works \cite{B1, B2, BEHP}. To
begin with, we give a summarized historic on BSDEs : the linear version
of  equation $(E^{(\xi,f)}) $ has appeared long time ago, both as
the equation for the adjoint process in stochastic control (see e.g.
\cite{Bi}), as well as the model behind the Black and Scholes
formula for the pricing and hedging of options in mathematical
finance, see e.g. \cite{BlSc, Me}. Since the paper \cite{PP1}, where
the existence and uniqueness of solutions have been established for
the equation $(E^{(\xi,f)})$ with a uniformly Lipschitz generator
$f$  and a square integrable  terminal data $\xi$, the theory of
BSDE has found further important applications and has become a
powerful tool in many fields such financial mathematics, optimal
control and stochastic game, non-linear PDEs ... etc. The collected
texts \cite{EM} give a useful introduction to the theory of BSDEs
and some of their applications. See also \cite{BEP, BM,  BPS, BL,
BP1, BH, EKPQ, HO, NOS, OT, P1, P2, R} and the references therein
for more discussions on BSDEs and their relations  with PDEs and
mathematical finances. Many authors have attempted to improve the
result of \cite{PP1} by weakening the Lipschitz continuity of the
coefficient $f$ (see e.g \cite {B1,B2, BEHP, Mao}) and/or the
$L^2$-integrability of the initial data $\xi$ (\cite {BDHPS, EKPQ}).
Another direction in the BSDEs theory has been  developed by
introducing the notion of \it{weak solutions}, i.e. a solution which
could be not adapted to the filtration generated by the driver
processes (see e.g. \cite {BMNO, BER, BE}). A step
forward has been done in the paper \cite {BER} where the
Meyer-Zheng topology has been used, to prove the existence of  weak
solutions for BSDEs with continuous generator.  More recently, the link between the
solution of BSDEs and the "$L_p-$\it{viscosity solution}" for
PDEs with discontinuous coefficients, has been established in \cite{BEP}.

The essential difficulty, to establish the existence of  strong
(i.e. $\mathcal{F}_t^W-$adapted) solutions to   BSDEs
with  local conditions on the
generator $f$, is due  to the fact that the
control variable
$Z$ is known implicitly, by Itô's martingale representation theorem as the
integrand
of a Brownian stochastic integral. Actually, we need more information on the
variable $Z$. Consequently, the usual localization
procedure (by stopping times) does not work.  On the other hand, the methods
used to study the existence and/or uniqueness of strong solutions to
one-dimensional BSDEs  are mainly based on comparison techniques   and therefore
do not work for multidimensional equations. We cite only a few articles in this
area (for instance \cite{BCEK, EH, FJ, K}) sending the reader to the references
therein  again because neither do
we deal with one-dimensional BSDEs nor use the results of these papers.  Note
however
 that, although we are focused in the multidimensional equations,
our  uniqueness result is new even in the one-dimensional case.
The first results which deal with the existence and
uniqueness as
well as the stability of strong solutions  for multidimensional BSDEs
{\it{with local
assumptions on the coefficient $f$}}
 have been  established in \cite{B1,B2,BEHP}.

 The present work constitute a natural development of  \cite{B1, B2,BEHP,BEH}.
 To begin, let $\xi$  be $p-$integrable with $p>1 $, $K$ be a positive
 constant,
 and
 consider the following example of  $d-$dimensional BSDE with logarithmic
 nonlinearity,
  \begin{equation}\label{EDSRlog}
  Y_t = \xi - \dint_t^T K Y_s \log \vert Y_s \vert ds - \dint_t^T Z_s dW_s
  \quad\quad
0\leq t\leq T
\end{equation}
 It is worth noting that the coefficient $f(y) := - K y \log \vert y \vert$ \, of  equation (\ref{EDSRlog}),  is not locally monotone and hence not locally
 Lipschitz. Moreover, its growth is big power than $y$.
 In our knowledge, when $\xi$ is $p-$integrable with $1 < p < 2 $, there is no
 results on multidimensional BSDEs which cover this interesting
 example.
To explain how the BSDE (\ref{EDSRlog})  follows naturally from
\cite{B1,B2}, consider the BSDE $(E^{(\xi,f)})$ with square
integrable $\xi$ and, assume for the simplicity that the generator
$f$ does not depend on the variable $z$. Let $f$ be $L_N-$locally
Lipschitz and with sublinear growth.  It has been established in
\cite{B1,B2} that if $L_N$ behaves as $\log N$, then the BSDE
$(E^{(\xi,f)})$ has a unique strong solution which is $L^2-$stable.
Now, if we drop the sublinear growth condition on $f$, then the
condition \ $L_N \sim \log N$ implies that $\vert f(y) \vert \leq
K(1 + \vert y \vert \log \vert y \vert)$ for some positive constant
$K$. Hence,  the following questions arise : could the BSDE with generator $
f(y) = - y  \log \vert y \vert$ has a strong solution ? If yes, what
happens about the uniqueness and the stability of solutions? These
questions are  positively solved  in this paper, as particular
examples.

The first main purpose of the present work consists to establish a result on the
existence and   uniqueness as well as the stability of strong solutions to
BSDE $(E^{(\xi,f)})$  which cover equation (\ref{EDSRlog}) as well as, other
interesting examples which are, in
our knowledge, not covered by the previous works.  For instance, we establish
the existence and uniqueness of
 a strong solution to BSDE $(E^{(\xi,f)})$ in the case where the terminal data
 $\xi$ is
merely $p$-integrable (with $p>1$) and the coefficient $f$ could be
neither locally monotone in $y$ nor locally Lipschitz in $z$.
Moreover,  $f$ can has a super-linear growth in its two variables
$y$ and $z$. For example, $f$ can take the form $f(y,z) =
-y\log\vert y\vert + g(y)(h(z)\sqrt{\vert \log\vert z\vert\vert})$
for some functions $g$ : $\R^d\longmapsto \R^d$ and $h$ :
$\R^d\times\R^r\longmapsto \R^d$. The assumptions which we impose on
$f$ are local in $y, z$ and also in $\omega $. This enables us to
cover certain BSDEs with stochastic monotone generators. Our
uniqueness result is  new even in the one-dimensional case.

The BSDEs with $p-$integrable terminal data $\xi$ (with $1<p<2$), have been
studied in  \cite{EKPQ} in the case where the coefficient $f$ is uniformly
Lipschitz in their two variables $(y,z)$, and in \cite{BDHPS} in the case
where  $f$ is uniformly Lipschitz in the variable $z$ and uniformly monotone
in the variable $y$. It should be noted that our result cover those of
\cite{BDHPS, EKPQ} with new proofs. Our method allows, for instance,
to treat simultaneously the existence and uniqueness as well as, the
$L^p-$stability of solutions by using the same computations.


   The  techniques  which is usually developed in BSDEs consist to applying
   It\^o's formula to the function $h(y) = \vert y\vert^2$ or
   $h(t,y) = \vert y\vert^2\exp(\alpha t)$ with $\alpha > 0$ in order to estimate
   the difference between two solutions by the difference between their respective
   data.  Such estimates are not possible in our situation since our assumptions
   on the generator are merely local. Moreover, due to the super-linear growth
   and the singularity of the generator, the  techniques used in \cite{B1,B2}
   can not be easily extended to our situation. Our proofs mainly consist to
   establish  a non standard a priori estimate  between two solutions by applying
   Itô's formula to an appropriate function. The existence (of solutions) is then
   deduced by using a suitable  approximation $( \xi_n, f_n)$ of $( \xi, f)$ and
   an appropriate
 localization procedure which is close to those given in \cite{B1,B2,BEHP}.
 However, in contrast to
\cite{BEHP}, we don't use the $L^2$-weak compactness of the approximating sequence
$(Y^n, Z^n)$. We directly show  that the sequence  $(Y^n, Z^n)$  strongly converges
in
some $L^{q}$ space ($1<q<2$) and, the limit satisfies the BSDE $(E^{(\xi,f)})$.
The uniqueness as well as the stability of solutions are then deduced by using
the same estimates. The results are first established  for a small time,  and next,
for an arbitrarily prescribed time  duration by using a continuation procedure.



 To deal with the PDEs part, we consider the  Markovian version of the BSDE
 (\ref{EDSRlog}) which is defined for $0 \leq  t\leq s\leq T$ by the system
 of SDE-BSDE,
\begin{equation}\label{E2}
\left\{
\begin{aligned}
&X^{t,x}_s  = x + \dint_t^s b(X^{t,x}_r) dr +  \dint_t^s
\sigma(X^{t,x}_r)dW_r, \\
&Y_s =H(X_t)-\int_s^t
K Y_r\log  \vert
 Y_r\vert dr-\int_s^t
Z_r dW_r
\end{aligned}
\right.
\end{equation}
 where \ $\sigma : \R^k\mapsto\R^{kr}$,\quad
$b : \R^k\mapsto\R^d$, \ \  $H : \R^k\mapsto\R^d$  are measurable
functions and $K$ is a real positive number.

\noindent The system of PDEs associated to the SDE-BSDE  (\ref{E2}) is then given
by
\begin{equation}\label{PDE1}
    \dfrac{\partial u(t,x)}{\partial t} + {L}u(t,x)- K u(t,x)\log  \vert
 u(t,x)\vert=0  \  ,  \ \ \ \ \ \ \ \
 u(T,x) = g(x)
\end{equation}
 where ${L}:=
\dfrac{1}{2}\displaystyle\sum_{i,j}(\sigma\sigma^*)_{ij}\partial2_{ij}+
\displaystyle\sum_i b_i\partial_i $ \ \ \ and \ \ \ $g$ is a given measurable
function.

\noindent
The logarithmic nonlinearities \ $Ku\log \vert u\vert$ \ [of the equation
(\ref{PDE1})] appear in some PDEs related to physics, see e.g.
\cite{BM1,BM2, CH, H,   PZ}. In the mathematical point of view, as indicated
in \cite{CH}, the nonlinear term $u\log \vert u \vert$ is not continuous on a
reasonable functions space. This induces a supplementary difficulty which makes
no efficient some standard arguments (local existence and global estimates)
to prove  existence of solutions. On the other hand, it should be noted that,
due to the degeneracy of the diffusion coefficient, the solutions
will not be smooth enough, and therefore the uniqueness is rather hard to
establish.

In the second part of this paper, we are concerned with the probabilistic
approach to Sobolev solutions of semilinear PDEs associated with the Markovian
version of BSDE $(E^{(\xi,f)})$. The links between  strong solutions of  BSDEs
and Sobolev solutions of semilinear PDEs were firstly established in \cite{BL}.
Similar result was established in \cite {BM}, for the relations between  Backward
Doubly SDEs (BDSDEs) and SPDEs. The common of these two papers is that the
nonlinear term $f$ is at least uniformly Lipschitz and with sub-linear growth.

  The second main purpose of this paper consists to establish a result on the
  existence and uniqueness of Sobolev solutions for the (possibly  degenerate)
  system of PDEs associated to BSDE $(E^{(\xi,f)})$. Our result cover equation
  (\ref{PDE1}) and many other examples. We develop a  method which allows to prove
  the uniqueness of the PDE by means of the uniqueness of its associated BSDE :
  we first prove the existence and uniqueness in the class of  solutions which are
  representable by  BSDEs, and next we show that any solution is
  unique. To do this, we first prove that 0 is the unique solution to the homogeneous
  PDE,
and next we use the BSDEs to establish an equivalence between the uniqueness for the non-homogeneous
semilinear PDE and the uniqueness for its associated homogeneous linear PDE. More precisely, denoting by  $\mathcal{L}$ the second
order parabolic
operator associated to a given $\R^d$-diffusion process, we prove that the
system of semilinear PDEs
\vskip0.2cm\hskip 1.5cm
$ \left\{
\begin{tabular}{ll}
 $\dfrac{\partial u(t,x)}{\partial t} + \mathcal{L}u(t,x)+f(t,x,u(t,x),
\nabla
 u(t,x))=0$,  \ \ \ \ \ \ \  $t\in ]0,T[$, $x\in \R^k$
 \\ $u(T,x) = g(x), \qquad$ $x\in\R^k$
\end{tabular} \right.$

\noindent has a unique solution {\it{if and
only if} $0$ is the unique solution of the linear system

\vskip0.2cm\hskip 1.5cm
$ \left\{
  \begin{tabular}{ll}
 $\dfrac{\partial u(t,x)}{\partial t} + \mathcal{L}u(t,x)=0$,   \ \ \ \ \ \ \
 $t\in
 ]0,T[$, $x\in \R^k$
 \\ $u(T,x) = 0, \qquad$ $x\in\R^k$
\end{tabular} \right.$

\vskip0.2cm\noindent This seems to be new in the  BSDEs framework.
Not also that, in order to prove the uniqueness of the above homogeneous linear
PDEs, a uniform gradient estimate for some possibly degenerate PDEs is established
by a probabilistic method, which is interesting itself.

\vskip 0.1cm We mention some others considerations which have
motivated the present work.

$\bullet$  The growth conditions on the nonlinearity constitute a critical case
in the sense that, for any $\varepsilon >0$, the solutions of the ordinary
differential equation $ X_t = x + \int_0^t X_s^{1+\varepsilon}ds$ \ explode
at a finite time.

$\bullet$   The logarithmic nonlinearities appear in some PDEs
arising in physics, see e.g. \cite{BM1,BM2, CH, H,  PZ, S}.
For
instance, in \cite{BM1}  the  construction  of nonlinear wave
quantum mechanics, based on Schrödinger-type equation, is with
nonlinearity $-k u \ln (|  u |^2)$. This nonlinearity is selected by
assuming the factorization of wave functions for composed systems.
Its most attractive features are : existence of the lower energy
bound. Moreover, it is the only one nonlinearity satisfying the
validity of Planck's relation  $E[\psi]=\hslash \ \psi$ for
stationary states $\psi$.

 $\bullet$
In terms of continuous-state branching processes, the logarithmic nonlinearity $u \log u$ corresponds to the Neveu branching mechanism.  This process was introduced by Neveu in \cite{N}, and  further studied in \cite{BeLeg, FS, FW}. For instance, the super-process with Neveu's branching mechanism constructed in \cite{FS} is related to the Cauchy problem,
\begin{equation}\label{PDEFleishman}
\left\{
\begin{aligned}
&\dfrac{\partial u}{\partial t} - \Delta u + u\log u =  0   \ \ \hbox{on} \ \ (0,\ \infty)\times \R^d  \\
&u(0^+) = \varphi >0
\end{aligned}
\right.
\end{equation}
Hence, our result can be seen as an alternative approach to the PDEs (\ref{PDEFleishman}), and cover the case where the diffusion part is possibly degenerate.

$\bullet$ Since the system of PDEs associated to the Markovian version of the  BSDE $(E^{(\xi,f)})$ can be degenerate, our result also covers certain systems of first order PDEs.

$\bullet$ Thanks to the possible degeneracy of the diffusion
coefficient, our proposition 4.2 cover for instance the PDE studied
in \cite{S} which arises in studying the motion of a particle acting
under a force perturbed by noise.

$\bullet$ The method, which we develop to study  the system of semilinear PDEs, is based on BSDEs and, both our results as well as their proofs  are new, particularly  the proof of the uniqueness.

$\bullet$ The BSDEs as well as the PDEs which we consider are interesting in themselves since the nonlinear part $f(t,y,z)$ can be neither locally monotone in $y$ nor locally Lipschitz in $z$.  Moreover, $f$ can be big power than $y$ and $z$, and therefore it is not uniformly continuous in $(y,z)$.

$\bullet$  It is worth noting that our condition on the coefficient $f$ is
 new even for the classical It\^o's forward SDEs. For instance,  we do not know whether or not the following  equation  (\ref{EDSxlogx}) possesses a pathwise unique solution.
 \begin{equation}\label{EDSxlogx} X_s = x + \dint_0^s X_r\log\vert X_r \vert dr +  \dint_0^s
X_r \sqrt{\vert \log\vert X_r\vert\vert} dW_r, \quad\quad\; 0\leq s\leq T
\end{equation}
It should be noted that the SDE (\ref{EDSxlogx}) is not covered by \cite{FZ}. We think that the method developed in the present paper may be used to solve this question. We are currently working on the SDEs (\ref{EDSxlogx}) since the stochastic flows of homeomorphisms defined by these type of SDEs seem be related to the construction of a metric in the Holder-Sobolev space $\mathcal{H}^{\frac{3}{2}}$, see \cite{M}.



\vskip 0.15cm The paper is organized as follows. In section 2, we present the assumptions and
the main result of the first part. We also give some examples. Section 3 is devoted to
the proof of the main result.  In section 4, we deal with PDEs: we  study the existence and uniqueness of a weak (Sobolev) $p-$integrable solution to systems of
degenerate semilinear PDEs whose nonlinearities are big power than  $u$ and $\nabla u$. We also
establish, in this section, an equivalence between the uniqueness for non-homogeneous semi-linear PDEs and the uniqueness for its associated homogeneous linear PDEs.


%
%
\section{\bf{Definition, assumptions, main result and examples.}}
%
%
 Throughout this paper, $p>1$ is an arbitrary fixed real number and all the
considered processes are $(\F_t)$-predictable.
\subsection{\bf{Definition.} }
  A solution of equation $(E^{(\xi,f)})$ is an
$(\F_t)$-adapted and $\R^{d+dr}$-valued process $(Y,Z)$ such that
$$\E\bigg(\sup_{t\leq T}|Y_t|^p +\big(\dint_0^T | Z_s | ^2 ds\big)^
{\frac{p}{2}} +\int_0^T|f(s,Y_s,Z_s)|ds \bigg) < +\infty $$ and satisfies
$(E^{(\xi,f)})$.
%
%
\subsection{\bf{Assumptions}}

\ni We consider the following assumptions on $(\xi, f)$:

\vskip0.5cm\noindent
There exist $ M\in\L^0(\Omega; \L^1([0, T]; \R_+))$,
 $ K \in \L^0(\Omega; \L^2([0, T]; \R_+))$
 and $\gamma\in ]0,\dfrac{1\wedge (p-1)}{2}[ $,  such that (with $\lambda_s:=2M_s + \dfrac{K_s^2}{2\gamma}$ ) we have,

\vskip0.3cm\noindent
\bf{(H.0)}\;
 \quad
  $ \E\mid\xi\mid^p\,e^{\frac{p}{2}\int_0^T \lambda_sds}
 <\infty$,
 \vskip 0.4cm \noindent\bf{(H.1)}\; $f
$ is continuous in $(y,z)$ for almost all $(t,\omega)$ \vskip0.4cm
\noindent
\bf{(H.2)}\;  There exist $\eta$ and $f^0 \in\L^{0}
(\Omega\times [0, T]; \R_
 +)$ satisfying

\par\hskip 0.2cm
 $\E\big(\dint_0^Te^{\int_0^s\lambda_rdr}\eta_s ds\big)^\frac{p}{2}<\infty$ , \ \ \ \
\quad $\E\big(\dint_0^Te^{\frac{1}{2}\int_0^s\lambda_rdr}f^0_s ds\big)^p<\infty$

 \noindent and such that :

 \vskip 0.2cm \ \ \ \ \ \ for every $t,y,z$, \ \  $ \langle y,f(t,y,z)\rangle \;\;\leq\;\;
\eta_t + f^0_t \vert y\vert + M_t\vert y\vert^2 + K_t\vert y\vert\vert z\vert $

\vskip0.4cm \noindent \bf{(H.3)}\;
 There exist $ \overline{\eta}\in\L^{q}(\Omega\times [0, T]; \R_
 +))$ (for some $q>1$) and $\alpha\in ]1,p[, \alpha'\in ]1,p\wedge 2[$
   such that:

   \vskip 0.2cm \ \ \ \ \ \ for every $t,y,z$, \ \  $\mid f(t,\omega, y, z)\mid \;\;\leq\;\;
\overline{\eta}_t +\mid y\mid^{\alpha}+\mid z\mid^{\alpha'}$.

\vskip0.4cm \noindent\bf{(H.4)}\;
 There exist $ v\in\L^{q'}(\Omega\times [0, T]; \R_
 +))$ (for some $q'>0$) and $K' \in\R_+$
 such that \\ for every
$N\in\N$  and every $  y,\; y'\; z,\; z' $  satisfying $ \mid y\mid,\; \mid
y'\mid,\; \mid z\mid, \;\mid z'\mid \leq  N $
\\ \\ $ \langle y-y^{\prime},
f(t,\omega,y,z)-f(t,\omega,y^{\prime},z')\rangle\1_{\{v_t(\omega)\leq N\}}$  $\leq K'\mid
y-y^{\prime}\mid^{2}log A_{N} + \sqrt{K'\log A_{N}}\mid y-y^{\prime}\mid\mid
z-z^{\prime}\mid {\hskip 7.5cm} +K' \dfrac{\log A_{N}}{ A_{N}}$\\ where $ A_N$ is a increasing
sequence and satisfies $A_N>1$, $\lim_{N\rightarrow\infty}{A_N}=\infty$ and
$A_N \leq N^{\mu}$ for some $\mu > 0$.

%
%

\subsection{\bf{The main result}}
\begin{theorem}\label{the21}
Assume that \bf{(H.0)-(H.4)} hold.   Then, $ (E^{(\xi, f)})$ has a unique solution $
(Y,Z)$ which satisfies,
\begin{align*}
 \E\sup_t\mid Y_t\mid^p & e^{\frac{p}{2}\int_0^t \lambda_sds} + \E\big[\dint_0^T
e^{\int_0^s \lambda_rdr}\mid Z_s\mid^2ds\big]^\frac{p}{2}  \\  \leq & C
\bigg\{\E\mid \xi\mid^pe^{\frac{p}{2}\int_0^T \lambda_sds} +\E\big(\dint_0^T
e^{\int_0^s \lambda_rdr}\eta_s ds\big)^{\frac{p}{2}}+ \E\big(\dint_0^T
e^{\frac{1}{2}\int_0^s \lambda_rdr}f^0_s ds\big)^{p} \bigg\}
\end{align*}
for some constant $C$ depending only on $p$ and $\gamma$.
\end{theorem}

We shall give some examples of BSDEs which satisfy the assumptions of Theorem
\ref{the21}. In our knowledge, these examples are not covered by the previous
works in multidimensional BSDEs.

\subsection{\bf{Examples.}}

{\bf{Example 1.}} Let $f(y):=-  y \log{\mid y \mid}$ then for all
$\xi\in\L^p(\F_T)$ the following BSDE has a unique solution $$ Y_t = \xi  -
\dint_t^T Y_s \log{\mid Y_s \mid}ds - \dint_t^T Z_s dW_s. $$ Indeed, $f$
satisfies \bf{(H.1)-(H.3)} since $ \langle y, f(y) \rangle \leq 1 $ and $\mid
f(y)\mid \leq 1+\dfrac{1}{\varepsilon}\mid y \mid^{1+\varepsilon} $ for all
$\varepsilon > 0.$ In order to verify \bf{(H.4)}, thanks to triangular
inequality, it is sufficient to treat separately the two cases: $0\leq \mid
y\mid,\mid y'\mid \leq \dfrac{1}{N}$ and $\dfrac{1}{N}\leq \mid y\mid,\mid
y'\mid \leq N$. \\ In the first case, since the map $x\mapsto -x\log x$
increases for $x\in]0,e^{-1}]$, we obtain for $N>e$
\begin{align*} |f(y)-f(y')|
&\leq \vert f(y)\vert + \vert f(y')\vert \\ &\leq 2\dfrac{\log{N}}{N}
\end{align*}
\\
In the second case, the finite increments theorem applied to $f$  shows that
\begin{align*}|f(y)-f(y')| \leq (1+\log{N})  \mid y-y'\mid.
\end{align*}
Hence \bf{(H.4)} is satisfied for every $N>e$ with $v_s=0$ and $A_N=N$.
\\\\ {\bf{Example 2.}} Let $ g(y):= y\log{\dfrac{\mid y\mid}{1+\mid y\mid}}$
and $h\in\mathcal{C}(\R^{dr};\R_+)\bigcap\mathcal{C}^1(\R^{dr}-\{0\};\R_+)$ be
such that

\vskip 0.2cm \hskip 1.5cm $h(z)$ = $ \left\{
\begin{tabular}{ll}
 $|z|\sqrt{-\log|z|}$ \hskip 1cm if $|z|< 1- \varepsilon_0$
\\ $|z|\sqrt{\log|z|}$ \hskip 1.2cm if $|z|> 1+ \varepsilon_0$
\end{tabular} \right.$

\vskip 0.2cm\noindent
 where $\varepsilon_0\in ]0,1[$.
Finally, we put $ f( y, z):= g(y)h(z)$.
Then for every $\xi\in\L^p(\F_T)$ the following BSDE has a unique solution
\begin{align*}
Y_t =  \xi  + \dint_t^T  f(Y_s,Z_s) ds - \dint_t^T Z_s dW_s.
\end{align*}
It is not difficult to see that $f$ satisfies (H1). We shall prove that $f$
satisfies (H2)-(H4).
\\

$(i)$ Since $g$ is continuous, $g(0)=0$ and $\vert g(y)\vert$ tends to $1$ as $\vert
y\vert$ tends to $\infty$, we deduce that $g$ is bounded. Moreover, $g$
satisfied $\langle y-y', g(y)-g(y')\rangle \leq 0$. Indeed, in one dimensional
case it is not difficult to show that $g$ is a decreasing function. Since,
$-\langle y, y'\rangle\log \frac{\vert y\vert}{1+\vert y\vert} \leq - \vert
y\vert\vert y'\vert\log \frac{\vert y\vert}{1+\vert y\vert}$ (because $\log
\frac{\vert y\vert}{1+\vert y\vert} \leq 0$), we can reduce the
multidimensional case to the one dimension case by developing the inner product
as follows,
\begin{align*} \langle y-y',
g(y)-g(y')\rangle & \leq \vert y\vert^{2} \log \dfrac{\vert y\vert}{1+\vert
y\vert} + \vert y'\vert^{2}\log \dfrac{\vert y'\vert}{1+\vert y'\vert} - \vert
y\vert\vert y'\vert(\log \dfrac{\vert y\vert}{1+\vert y\vert} + \log
\dfrac{\vert y'\vert}{1+\vert y'\vert})
\\ &
= (\vert y\vert -\vert y'\vert ) (\vert y\vert\ \log \dfrac{\vert
y\vert}{1+\vert y\vert} - \vert y'\vert\ \log \dfrac{\vert y'\vert}{1+\vert
y'\vert})  \\& = \langle \vert y\vert -\vert y'\vert , g(\vert y\vert )-g(\vert
y'\vert )\rangle \\ & \leq 0
\end{align*}

 $(ii)$  The function $h(z)$ satisfies for all $
 \varepsilon >0$
 \begin{align*}0 \leq h(z)\leq M + \dfrac{1}{\sqrt{2\varepsilon}} \mid
z\mid^{1+\varepsilon}, \qquad \hbox{where}\quad M = \sup_{\mid z\mid \leq
1+\varepsilon_0}{\mid h(z) \mid}
\end{align*}
The last inequality follows since
 $\sqrt{2\varepsilon\log\vert z\vert}$ =
 $\sqrt{\log\vert z\vert^{2\varepsilon}} \leq \vert z\vert^{\varepsilon}$
 for each $\varepsilon>0$ and $\vert z\vert>1$.
 (H3)  follows now
 directly from the previous observations (i) and (ii). (H2) is satisfied since
 $\langle y , f(y,z)\rangle  = \langle y , g(y)
 \rangle  h(z) \leq 0 $.
To verify \bf{(H.4)} it is enough to show that for every $ z,z'$ such that $
\mid z\mid, \mid z'\mid \leq N$ $$ \mid h(z)- h(z')\mid \leq
c\left(\sqrt{\log{N}}\mid z-z'\mid + \dfrac{\log{N}}{N}\right)$$ for $N$ large enough and some
positive constant $c$. This can be proved by considering
separately the following five cases,  $ 0 \leq \mid z\mid, \mid z'\mid \leq
\dfrac{1}{N}$,\quad $ \dfrac{1}{N} \leq \mid z\mid, \mid z'\mid \leq
1-\varepsilon_0 $,\quad $ 1-\varepsilon_0 \leq \mid z\mid, \mid z'\mid \leq
1+\varepsilon_0 $ and   $ 1+\varepsilon_0 \leq \mid z\mid, \mid z'\mid \leq N
$.
\\\\
 In the first case (i.e.
$ 0 \leq \mid z\mid, \mid z'\mid \leq \frac{1}{N}$), since the map $x\mapsto
x\sqrt{-\log x}$ increases for $x\in [ 0,\frac{1}{\sqrt e}]$, we obtain  \  \quad $ |h(z)-h(z')| \leq \vert h(z)\vert + \vert h(z')\vert \leq
2\dfrac{1}{N}\sqrt{-\log{\dfrac{1}{N}}} \leq 2\dfrac{1}{N}\log N $ \ for
$N>\sqrt e$.
\\The other cases can be
proved by using the finite increments theorem.
\\\\


\noindent {\bf{Example 3.}} Let $(X_t)_{t\leq T}$ be an $(\F_t)-$adapted and
$\R^k-$valued process satisfying the forward stochastic differential
equation $$ X_t = X_0 + \dint_0^t b(s,X_s)ds + \dint_0^t \sigma (s, X_s) dW_s
$$ where $X_0\in\R^k$ and $ \sigma, b \; :\, [0, T]\times\R^k\rightarrow
\R^{kr}\times\R^k$ are measurable functions such that $ \|\sigma(s,x)\| \leq c
$ and $ |b(s, x)|\leq c(1+|x|)$, for some constant $c$.


It is known from the forward SDE's theory that  there exist $ \kappa
>0$ and $C
> 0$ depending only on $c, T, k$ such that $$ \E\exp{(\kappa
\sup_{t\leq T}{\mid X_t\mid^2})} \leq C \exp(C\mid X_0\mid^2).$$

\noindent Consider the BSDE

$$ Y_t = g(X_T) + \dint_t^T \mid X_s\mid^{\overline{q}} Y_s - Y_s\log{\mid
Y_s\mid} ds - \dint_t^T Z_s dW_s.
$$
where $\overline{q}\in ]0,2[$ and $g$ is a measurable function satisfying $\mid
g(x)\mid \leq c\exp{c\mid x\mid^{\overline{q}'}}$, for some constants $c>0,
\;\overline{q}'\in [0,2[$.\\ The previous BSDE has a unique solution $(Y, Z)$
 which satisfies: for every $p>1$ there
exists a positive constant $C$ such that $$ \E\sup_t\mid Y_t\mid^p +
\E\left[\dint_0^T \mid Z_s\mid^2ds\right]^\frac{p}{2} \leq C \exp{(C\mid
X_0\mid^{2})}.$$

\noindent Indeed, one can show that\\ $i)$ $ \langle y, f(t,y)\rangle \leq 1 +
\mid X_t\mid^{\overline{q}} \mid y\mid^2$
\\ $ii)$ Using Young inequality we obtain, for every $\epsilon
>0 $ there is a constant  $c_{\epsilon} >0 $ such,
that $$ \mid f(t,y)\mid \leq c_{\epsilon} (1+\mid
X_t\mid^{\overline{q}c_{\epsilon}} +\mid y\mid^{1+\epsilon}) $$
\\
$iii)$ $ f $ satisfies assumption \bf{(H.4)} with $v_s = \exp{\mid
X_s\mid^{\overline{q}}}$ and $A_N=N$.


\vskip 0.6cm The following example shows that our assumptions enable to treat
BSDEs with stochastic monotone coefficient

\vskip 0.2cm\noindent\bf{Example 4.} Let $(\xi, f)$ satisfying \bf{(H.0)-(H.3)}
and
\\ \bf{(${H'}$.4)}\;  $\left\{
\begin{tabular}{ll}
 There are a positive process $C$ satisfying $ \E\dint_0^T
 e^{q'C_s}ds<\infty$ (for some $q'>0$)
 and $ K'\in\R_+ $\\ such that: \\
\\ $ \langle y-y^{\prime},
f(t,\omega,y,z)-f(t,\omega,y^{\prime},z')\rangle \leq  K'\mid
y-y^{\prime}\mid^{2}\left\{C_t(\omega) + \mid\log{\mid y-y'\mid}\mid \right\} $
\\ \\ \hspace{55mm} $ + K'\mid y-y^{\prime}\mid\mid z-z^{\prime}\mid \sqrt{
C_t(\omega) + \mid\log{\mid z-z'\mid}\mid}.$
\end{tabular} \right.$ \vskip0.5cm

 \noindent In particular we have for all $z,z'$ $$
|f(t,\omega,y,z)-f(t,\omega,y,z')| \leq  K'\mid z-z^{\prime}\mid \sqrt{
C_t(\omega) + \mid\log{\mid z-z'\mid}\mid}.$$ Therefore,  the following BSDE has a
unique solution $$ Y_t = \xi  + \dint_t^T f(s, Y_s, Z_s)ds - \dint_t^T Z_s
dW_s. $$
To check \bf{(H.4)}, it is enough to show that for some constant
$c$ we have

 \ni $ \langle y-y', f(t,y,z)-f(t,y',z)\rangle \leq c \log{N}\left(
\mid y-y'\mid^2 + \dfrac{1}{N} \right)$ \\ $
 | f(t,y,z)-f(t,y,z')|\leq c \sqrt{\log{N}}\left(
\mid z-z'\mid + \dfrac{1}{N} \right)$

\noindent whenever $v_s:= e^{C_s}\leq N$ and $\mid y\mid$ , $\mid y'\mid$,
$\mid z\mid$, $\mid z' \mid$ $\leq N$.

\noindent These two inequalities can be respectively proved by considering the following cases $$
\mid y-y'\mid \leq \dfrac{1}{2N}, \;\; \;\;\dfrac{1}{2N}\leq \mid y-y'\mid \leq
2N. $$ and $$ \mid z-z'\mid \leq \dfrac{1}{2N}, \;\; \;\;\dfrac{1}{2N}\leq \mid
z-z'\mid \leq 2N. $$


\vskip 0.6cm\noindent\bf{Example 5.} Let $(X_t)_{t\leq T}$ and $\xi$ be as in
example 3, let $F(t, x, y, z)$ be such that \\ $i)$ $F(t,x,.)$ is continuous
\\ $ii)$ $ |F(t, x, y, 0)|\leq C\exp(C\mid x\mid^{\hat{q}})+\mid
y\mid^\alpha$, for some $ \hat{q},\alpha\in ]0,2[$ and $C>0,$
\\ $iii)$ $ \langle F(t, x, y, z)-F(t, x, y', z'),y-y'\rangle \leq
K'\mid y-y'\mid^2+K'\mid y-y'\mid\mid z-z'\mid$.
\\ Let $\overline{q},\overline{q}', \overline{q}"\geq 0$ such that $
\overline{q}+\overline{q}"<2 $ and $ \overline{q}'+\overline{q}"<1$, the
following BSDE has a unique solution $$ Y_t = \xi  + \dint_t^T
|X_s|^{\overline{q}"}F(s, X_s,|X_s|^{\overline{q}} Y_s, |X_s|^{\overline{q}'}
Z_s)ds - \dint_t^T Z_s dW_s. $$

\section{\bf{Proof of Theorem \ref{the21}} }
We first give some a priori estimates from which we derive a stability result
for BSDEs and next we use  a suitable approximation
of $(\xi, f)$ to complete the proof. The difficulty comes from the fact
that the generator $f$ can be neither locally monotone in the variable $y$ nor locally Lipschitz
in the variable $z$ and moreover, it also may have a super-linear growth in its
two variables $y$ and $z$.


\subsection{\bf{Estimates for the solutions of equation $(E^{(\xi,f)})$.}}
In the first step, we give estimates for the processes $Y$ and $Z$.


\begin{proposition}\label{pro41}   Let $\Lambda_t:= | Y_{t}|
^{2}e_t+2\displaystyle\int_0^t e_s\eta_s ds +(\int_0^t
e_s^\frac{1}{2}f^0_s ds)^{2}$ and $e_t:=\exp{\int_0^t\lambda_s}ds $.

\noindent Assume that \bf{(H.2)} hold and  \ $\displaystyle\E(\sup_{0\leq\ s \leq T} \left| Y_{t}\right|
^{p}e_t^\frac{p}{2})<\infty$.

\noindent Then, there exists a positive
constant $ C^{(p,\gamma)} $ such that
\begin{align*}
\displaystyle\E\sup_{0\leq\ s \leq T}\Lambda_s^\frac{p}{2} & +
\displaystyle\E\left(\displaystyle\dint_0^T e_s\vert
Z_s\vert^2ds\right)^{\frac{p}{2}}
 \leq  C^{(p,\gamma)}
\E\Lambda_T^\frac{p}{2}.
\end{align*}
\end{proposition}


To prove this proposition we need some lemmas.
\begin{lemma}\label{lem40}
For every $\varepsilon>0$, every $\beta>1$ and every positive functions $h$ and
$g$ we have $$\dint_t^T(h(s))^{\frac{\beta -1}{2}}g(s)ds\leq
\varepsilon\sup_{t\leq s \leq T}\mid
h(s)\mid^{\frac\beta2}+\varepsilon^{1-\beta
}\left(\dint_t^Tg(s)ds\right)^{\beta}.$$
\end{lemma}


\noindent\bf{Proof.} Let $\varepsilon>0$ and $\beta>1$. Using Young's inequality we get for
every $\delta$ and $\delta'$ such that $\frac{1}{\delta}+\frac{1}{\delta'}=1$
\begin{align*}
&\dint_{t}^{T}(h(s))^{\frac{\beta -1}{2}}g(s)ds\leq
\frac{1}{\delta}\varepsilon^{\frac{(\beta -1)\delta}{\beta}}\sup_{t\leq s \leq
T}\mid h(s)\mid^{\frac{(\beta -1)\delta}{2}}+\frac{\varepsilon^{\frac{(1-\beta
)\delta'}{\beta}}}{\delta'}\big(\dint_t^Tg(s)ds\big)^{\delta'}
\end{align*}
 We now choose $\delta =\frac{\beta}{\beta -1}$ and use
 the fact that $\delta, \delta' > 1$.\eop


\begin{lemma}\label{lem41}
If  \bf{(H.2)} holds then for  every  $\beta >1+2\gamma$  there exist positive
constants $C^{(\beta,\gamma)}_1  , C^{(\beta,\gamma)}_2 $ such that for  every
$\varepsilon
>0$, every stopping time $\tau \leq T$ and every $ t\leq \tau$
\begin{align*}
& \Lambda_t^\frac{\beta}{2}+\displaystyle\dint_t^\tau
\Lambda_s^\frac{\beta-2}{2}e_s\vert Z_s\vert^2ds
  \;\;\leq\;\;
\varepsilon\sup_{t\leq s \leq \tau}\Lambda_s^\frac{\beta}{2}+\varepsilon^{(1-
\beta) }C^{(\beta,\gamma)}_1 \Lambda_\tau^\frac{\beta}{2}
 -C^{(\beta,\gamma)}_2 \displaystyle\dint_t^\tau \Lambda_s^{\frac{\beta}{2}-1}e_s
\langle Y_s,Z_s dW_s \rangle.
\end{align*}
\end{lemma}


\vskip0.2cm \noindent\bf{Proof.} Without loss of generality, we assume that $ \eta $ and $
f^0 $ are strictly positives.\\ It follows by using It\^o's formula that for
every $t\in[0,\tau]$,
\begin{align*}
\left| Y_{t}\right| ^{2}e_t+\dint_{t}^{\tau}\left| Y_{s}\right|^{2}\lambda_se_s
ds
 &= e_{\tau}\left| Y_{\tau}\right|
^{2}+2\dint_{t}^{\tau}e_s\langle Y_{s},\,f(s,Y_{s},Z_{s})\rangle ds
 -\dint_{t}^{\tau}e_s\mid Z_{s}\mid^2
ds\\&\quad-2\dint_{t}^{\tau}e_s\langle Y_{s},\;
 Z_{s}dW_{s}\rangle.
\end{align*}
Again It\^o's formula, applied to the process $\Lambda$, shows that
\begin{align*}
\Lambda_{t}^{\beta\over 2} &+\beta\dint_{t}^{\tau}
\Lambda_{s}^{\frac{\beta}{2}-1}\left(\frac12\left| Y_{s}\right|^{2}\lambda_se_s
+ e_s\eta_s + f^0_se_s^{\frac{1}{2}} \left[\dint_0^s f^0_r e_r^{\frac{1}{2}}
dr\right]\right) ds\\& = \Lambda_\tau^{\beta\over 2}
+\beta\dint_{t}^{\tau}\Lambda_{s}^{\frac{\beta}{2}-1} \langle e_sY_{s}, \quad
f(s,Y_{s},Z_{s})\rangle ds
 -\frac{\beta}{2}\dint_{t}^{\tau}\Lambda_{s}^{\frac{\beta}{2}-1}\left|
Z_{s}\right|^{2}e_s ds
\\ &\quad-\beta\dint_{t}^{\tau}e_s\Lambda_{s}^{\frac{\beta}{2}-1}\langle
Y_{s}, \quad  Z_{s}dW_{s}\rangle
-\beta(\frac{\beta}{2}-1)\dint_{t}^{\tau}e_s^2\Lambda_{s}^{\frac{\beta}{2}-2}
\sum_{j=1}^{r} \left(\sum_{i=1}^{d}Y_{s}^i Z_{s}^{i,j}\right)^{2}ds
\end{align*}
Observe that $\displaystyle \sum_{j=1}^{r}\left(\sum_{i=1}^{d}Y_{s}^i
Z_{s}^{i,j})\right)^{2} \leq \vert Y_s\vert^2\vert Z_s\vert^2 \leq
e^{-1}_s\Lambda_{s}\left| Z_{s}\right|^{2} $ then use the assumption \bf{(H.2)}
to get
\begin{align*}
&\Lambda_{t}^{\beta\over 2} + \frac{\beta}{2}
(1-2\gamma-(2-\beta)^+)\dint_{t}^{\tau}\Lambda_{s}^{\frac{\beta}{2}-1}e_s\left|
Z_{s}\right|^{2} ds \\& \leq \;\; \Lambda_\tau^{\beta\over 2}
 +\beta\dint_{t}^{\tau}\Lambda_{s}^{\frac{\beta}{2}-\frac{1}{2}}
 f^0_se_s^{\frac{1}{2}}  ds
-\beta\dint_{t}^{\tau}\Lambda_{s}^{\frac{\beta}{2}-1}\langle e_sY_{s}, \quad
Z_{s}dW_{s}\rangle.
\end{align*}
It follows from Lemma \ref{lem40} with $ g(s) = f^0_se_s^{\frac{1}{2}}$, since
$ \left( \dint_{t}^{\tau}
 f^0_se_s^{\frac{1}{2}} ds\right)^{\beta}  \leq \;\Lambda_{\tau}^
 {\frac{\beta}{2}},$ that for  every
$ \varepsilon
>0 $
$$ \dint_{t}^{\tau}\Lambda_{s}^{\frac{\beta}{2}-\frac{1}{2}}
 f^0_se_s^{\frac{1}{2}}  ds
 \leq \;\; \varepsilon \sup_{t\leq s\leq \tau}
 \Lambda_{s}^{\frac{\beta}{2}}+\varepsilon^{1-\beta}
 \Lambda_{\tau}^{\frac{\beta}{2}}
$$
Since $\beta > 1+2\gamma$ implies that $1-2\gamma - (2-\beta)^+ > 0$, Lemma
\ref{lem41} is proved.\eop


\begin{lemma}\label{lem42} Let \bf{(H2)} be satisfied and assume that
$\E(\sup_{0\leq s \leq T} \left| Y_{t}\right|
^{p}e_t^\frac{p}{2})<\infty$.

\noindent Then,
\\
$ 1)$There exists a positive constant $ C^{(p,\gamma)} $ such that for every  $
\varepsilon
>0$, we have
\begin{align*}
\displaystyle\E\displaystyle\int_0^T & \Lambda_s^\frac{p-2}{2}e_s\vert
Z_s\vert^2ds
 \;\;\leq\;\;
\varepsilon\E(\sup_{0\leq s \leq
T}\Lambda_s^\frac{p}{2}) \ + \ \varepsilon^{(1-p)}C_1^{(p,\gamma)}
\E(\Lambda_T^\frac{p}{2}).
\end{align*}
$ 2)$ There exists a positive constant $ C^{(p,\gamma)} $ such that
   $$\displaystyle\E\big(\displaystyle\int_0^T
e_s\vert Z_s\vert^2ds\big)^{\frac{p}{2}} \;\;\leq\;\; C^{(p,\gamma)}
\displaystyle\E(\sup_{0\leq s \leq T}\Lambda_s^\frac{p}{2}).$$
\end{lemma}


\noindent\bf{Proof.} The first assertion follows by a standard martingale localization
procedure. To prove the second assertion, we successively use Lemma \ref{lem41}
(with $\varepsilon =1$ and $\beta=2$), the Burkholder-Davis-Gundy inequality,
the fact that $e_s \vert Y_s\vert^2 \leq \Lambda_s $ and Young's inequality to
obtain
\begin{align*}
 \E\big(\displaystyle\dint_0^T e_s\vert
Z_s\vert^2ds\big)^{\frac{p}{2}}
 & \leq
C^{(p,\gamma)}_1\E\big(\sup_{0\leq s \leq T}\Lambda_s^\frac{p}{2}\big)
 +C^{(p,\gamma)}_2 \E\big (\vert \displaystyle\dint_t^T
 e_s \langle
Y_s,Z_s dW_s \rangle\vert^{\frac{p}{2}}\big)\\ & \leq
C^{(p,\gamma)}_1\E\big(\sup_{0\leq s \leq T}\Lambda_s^\frac{p}{2}\big)
 + C^{(p,\gamma)}_2 \E\big (\vert\displaystyle\dint_0^T
  e^2_s\vert
Y_s\vert^2\vert Z_s\vert^2 ds\vert^{\frac{p}{4}}\big)\\ &\leq
C^{(p,\gamma)}_1\E\big(\sup_{0\leq s \leq T}\Lambda_s^\frac{p}{2}\big)
 + C^{(p,\gamma)}_2  \E\big (\vert\displaystyle\dint_0^T
 \Lambda_se_s
 \vert Z_s\vert^2 ds\vert^{\frac{p}{4}}\big)\\ &\leq
 C^{(p,\gamma)}_1\E\big(\sup_{0\leq s \leq
T}\Lambda_s^\frac{p}{2}\big)
 + C^{(p,\gamma)}_2  \E\big [(\sup_{0\leq s \leq
T}\Lambda_s^{\frac{p}{4}})(\displaystyle\dint_0^T
 e_s
 \vert Z_s\vert^2 ds)^{\frac{p}{4}}\big]\\ &\leq
 \big [C^{(p,\gamma)}_1+
 + 2(C^{(p,\gamma)}_2 )^2\big]\E(\sup_{0\leq s \leq
T}\Lambda_s^{\frac{p}{2}}) + \frac{1}{2}\E\big[(\displaystyle\dint_0^T
 e_s
 \vert Z_s\vert^2 ds)^{\frac{p}{2}}\big]\\ &\leq
 [2C^{(p,\gamma)}_1+4(C^{(p,\gamma)}_2 )^2]\;\;\;\E(\sup_{0\leq s \leq
T}\Lambda_s^\frac{p}{2}).
\end{align*}
 \ni Lemma \ref{lem42}
 is proved.
\eop


\begin{lemma}\label{lem43}   Let the assumptions of Lemma 3.3 be satisfied.
Then, there exists a constant $ C^{(p,\gamma)} $ such that
$$\displaystyle\E(\sup_{0\leq s \leq T}\Lambda_s^\frac{p}{2}) \;\;\leq\;\;
C^{(p,\gamma)}\displaystyle\E(\Lambda_T^\frac{p}{2}) .$$
\end{lemma}


\noindent\bf{Proof.} Lemma \ref{lem41} and the Burkholder-Davis-Gundy inequality show that
there exists a universal positive constant $ c $ such that for every
$\varepsilon>0$ and  $t\leq T $
\begin{align*}
 \E\sup_{0\leq s \leq T}\Lambda_s^\frac{p}{2}\;\;\leq\;\; &
\varepsilon\E(\sup_{0\leq s \leq
T}\Lambda_s^\frac{p}{2})+\varepsilon^{(1-p)}C^{(p,\gamma)}_1
\E(\Lambda_T^\frac{p}{2})
\\
& +
 c C^{(p,\gamma)}_2 \E\big(\displaystyle\dint_0^T \Lambda_s^{p-2}
(|Y_s|^2e_s)e_s|Z_s|^2 ds \big)^{\frac{1}{2}}.
\end{align*}
Young's inequality gives, for  every  $\varepsilon'>0$,
\begin{align*}
 \E(\sup_{0\leq s \leq T}\Lambda_s^\frac{p}{2})\;\;\leq\;\; &
\varepsilon\E(\sup_{0\leq s \leq
T}\Lambda_s^\frac{p}{2}) \ + \ \varepsilon^{(1-p)}C^{(p,\gamma)}_1
\E(\Lambda_T^\frac{p}{2})
\\
& +\varepsilon'\E(\sup_{0\leq s \leq T}\Lambda_t^\frac{p}{2}) \ + \ \dfrac{\left[c
C_2^{(p,\gamma)}\right]^2 }{\varepsilon'}
 \E\displaystyle\dint_0^T \Lambda_s^\frac{p-2}{2}e_s\vert Z_s\vert^2ds.
\end{align*}
Applying Lemma \ref{lem42}, we get for every $\varepsilon">0$
\begin{align*}
 \E(\sup_{0\leq s \leq T}\Lambda_t^\frac{p}{2})\;\;\leq\;\; &
(\varepsilon+\varepsilon'+\dfrac{\left[c C^{(p,\gamma)}_2\right]^2\varepsilon"
}{\varepsilon'})\E(\sup_{0\leq s \leq T}\Lambda_s^\frac{p}{2})
\\
& + (\varepsilon^{(1-p)}C^{(p,\gamma)}_1 +\dfrac{\left[c
C^{(p,\gamma)}_2\right]^2
C^{(p,\gamma)}_1(\varepsilon")^{(1-p)}}{\varepsilon'}) \E(\Lambda_T^\frac{p}{2}).
\end{align*}
A suitable choice of $ \varepsilon, \varepsilon', \varepsilon" $ allows to conclude the
proof. \eop


\vskip 0.2cm\noindent\bf{Proof of Proposition \ref{pro41}.} It follows from Lemma 3.3 and
Lemma 3.4. \eop


\begin{proposition}\label{pro42} If \bf{(H.3)}
 holds  then,
\begin{align*}
&  \E  \dint_0^T   | f(s,Y_s,Z_s)|^{\hat{\beta}} ds  \; \leq \;
9^{p+q}(1+T)\big[1+\E\dint_0^T
 \overline{\eta}_s^{q} ds+ \E\sup_{0\leq s\leq T}\vert Y_s\vert^{p} +
\E(\dint_0^T\vert Z_s\vert^2 ds)^{\frac{p}{2}}\big]
\end{align*}
where $ \hat{\beta} :=
\dfrac{2}{\alpha'}\wedge\dfrac{p}{\alpha}\wedge\dfrac{p}{\alpha'}\wedge q$.
\end{proposition}


\noindent\bf{Proof.} We successively use Assumption (H.3), Young's inequality
and H\" older's inequality to show that
\begin{align*}
 \E  \dint_0^T   | f(s,Y_s,Z_s)|^{\hat{\beta}} ds
 &\leq
\E\dint_0^T (\overline{\eta}_s + \vert Y_s\vert^\alpha + \vert
Z_s\vert^{\alpha'})^{\hat{\beta}} ds
\\ &\leq 3^{\hat{\beta}}\E\dint_0^T
(\overline{\eta}_s^{\hat{\beta}} + \vert Y_s\vert^{\alpha\hat{\beta}} + \vert
Z_s\vert^{\alpha'\hat{\beta}}) ds
\\ &\leq
3^{\hat{\beta}}\E\dint_0^T ((1 + \overline{\eta}_s)^{\hat{\beta}} + (1 + \vert
Y_s\vert)^{\alpha\hat{\beta}} + (1 + \vert Z_s\vert)^{\alpha'\hat{\beta}}) ds\\
&\leq 3^{\hat{\beta}}\E\dint_0^T ((1 + \overline{\eta}_s)^{q} + (1 + \vert
Y_s\vert)^{p} + (1 + \vert Z_s\vert)^{p\wedge2}) ds
\\ &\leq 3^{\hat{\beta}}3^{p+q}\E\dint_0^T
(1 + \overline{\eta}_s^{q} + \vert Y_s\vert^{p} + \vert Z_s\vert^{p\wedge2})
ds\\ &\leq 3^{\hat{\beta}}3^{p+q}\big[T+\E\dint_0^T  \overline{\eta}_s^{q} ds +
T\E\sup_{0\leq s\leq T}\vert Y_s\vert^{p}
+T^{\frac{2-(p\wedge2)}{2}}\E(\dint_0^T\vert Z_s\vert^2 ds)^{\frac{p}{2}}\big]
\\
\\ &\leq 9^{p+q}(1+T)\big[1+\E\dint_0^T
 \overline{\eta}_s^{q} ds+ \E\sup_{0\leq s\leq T}\vert Y_s\vert^{p} +
\E(\dint_0^T\vert Z_s\vert^2 ds)^{\frac{p}{2}}\big].
\end{align*}
Proposition \ref{pro42} is proved. \eop


 \subsection{\bf{Estimate of the difference between two solutions.}}

 The next proposition gives an estimate which is a key tool in
the proofs.


\begin{lemma} \label{lem43}
Let $(\xi^i, f_i)_{i=1,2}$ satisfy \bf{(H.3)} (with the same $\overline{\eta},
\alpha$ and $\alpha'$) and let $(Y^i, Z^i)$ be a solution of $(E^{(\xi^i,
f_i)})$. Then, there exist $\beta = \beta(p,q,\alpha,\alpha')\in ]1,p\wedge
2[$, $r= r(p,q,\alpha,\alpha',K',\mu,q')>0$ and $a =
a(p,q,\alpha,\alpha',K',\mu,q') > 0$ such that for every $ u \in [0, T], u' \in
[u,T\wedge (u+r)]$, $N>0$ and every function $ f $ satisfying $\bf{(H.4)}$
\begin{align*}
& \E(\sup_{u\leq t\leq u'} |Y^1_t-Y^2_t|^{\beta}) \ + \ \E
\dint_{u}^{u'}\dfrac{|Z_{s}^{1}-Z_{s}^{2}|^2}
{\left(1+|Y^1_s-Y^2_s|^2\right)^{1-\frac{\beta}{2}}}
 ds
\\&
 \leq  \;\; NA_N^{1+\frac{\beta}{2}}\left[ \E
(|Y^1_{u'}-Y^2_{u'}|^{\beta}) + \E\dint_{0}^{T}\rho_{N} (f_1-f)_s+\rho_{N}
(f_2-f)_s  ds \right]
\\ & \ \ \ \ \ \ + \ \dfrac{1}{A_N^a}\left[1+\Theta_p^1+\Theta_p^2+\E\dint_0^T
\overline{\eta}_s^qds+ \E\dint_0^T v_s^{q'}ds\right].
\end{align*}
where $$ \rho_N(f_i-f)(t, \omega):= \sup_{\vert y\vert, \vert
z\vert\;\;\leq\;\; N}\vert f(t,\omega,y,z)-f_i(t,\omega,y,z)\vert$$ and $$
\Theta^i_p \ := \ \E(\sup_t|Y^i_t|^p) \ + \ \E\left(\dint_{0}^{T}|Z^i_{s}|^{2} ds
\right)^{\frac{p}{2}}.$$
\end{lemma}


\noindent\bf{Proof.} Let $q$ be the number defined in assumption \bf {(H3)} and $q',K', \mu$
those defined in assumption \bf {(H4)}. Let $\overline{\gamma}>0$ be such that
$ 1+2\overline{\gamma}<\hat{\beta}:= \dfrac{2}{\alpha'}\wedge
\dfrac{p}{\alpha}\wedge\dfrac{p}{\alpha'}\wedge q $ and set $K":=
K'+\dfrac{K'}{4\overline{\gamma}}$. Let $\beta \in ]1+2\overline{\gamma},
\hat{\beta}[$ and $ \nu \in ]0, (1-\frac{\beta}{\hat{\beta}})(1 \wedge q')[$.
Let $ r\in ]0,\dfrac{\nu}{\mu\hat\beta K"}\wedge\dfrac{1}{2K"}\wedge 1[$.
\\
\ni For \;$N\in \N $,   we set $$  \overline{e}_t:= (A_N)^{2 K"(t-u)}\quad
\hbox{and} \quad\Delta_{t}:=\{\left| Y_{t}^1-Y_{t}^2\right| ^{2}+
(A_N)^{-1}\}\overline{e}_t.$$
\\
%
%
 Using It\^{o}'s formula, we show that for every stopping time $\tau\in [u,
u']$ and every $ t \in [u,\tau]$
\begin{equation}\label{equ1}
 \begin{array}{lll}
\Delta_{t}^{\beta\over 2}&+&2\log (A_N) K"\dint_{t}^{\tau}\overline{e}_s
\Delta_{s}^{\frac{\beta}{2}}ds
+\frac{\beta}{2}\dint_{t}^{\tau}\overline{e}_s\Delta_{s}^{\frac{\beta}{2}-1}\left|
Z_{s}^{1}-Z_{s}^{2}\right|^{2}ds\\&&\\&& = \Delta_\tau^{\beta\over
2}-\beta\dint_{t}^{\tau}\overline{e}_s\Delta_{s}^{\frac{\beta}{2}-1}\langle
Y_{s}^{1}-Y_{s}^{2},\left( Z_{s}^{1}-Z_{s}^{2}\right)dW_{s}\rangle\\
&&\\&&\quad +\beta\dint_{t}^{\tau}\overline{e}_s\Delta_{s}^{\frac{\beta}{2}-1}
\langle Y_{s}^{1}-Y_{s}^{2},
f_{1}(s,Y_{s}^{1},Z_{s}^{1})-f_{2}(s,Y_{s}^{2},Z_{s}^{2})\rangle ds
\\ &&\\&&\quad
-\beta(\frac{\beta}{2}-1)\dint_{t}^{\tau}\overline{e}^2_s\Delta_{s}^{\frac{\beta}{2}-2}
\sum_{j=1}^{r}
\left(\sum_{i=1}^{d}(Y_{i,s}^{1}-Y_{i,s}^{2})(Z_{i,j,s}^{1}-Z_{i,j,s}^{2})\right)^{2}ds\\&&\\&&
=\Delta_\tau^{\beta\over
2}-\beta\dint_{t}^{\tau}\overline{e}_s\Delta_{s}^{\frac{\beta}{2}-1}\langle
Y_{s}^{1}-Y_{s}^{2},\left( Z_{s}^{1}-Z_{s}^{2}\right)dW_{s}\rangle+\beta I_1
-\beta(\frac{\beta}{2}-1)I_2,
\end{array}
\end{equation}
where
\[
I_1:=\dint_{t}^{\tau}\overline{e}_s\Delta_{s}^{\frac{\beta}{2}-1} \langle
Y_{s}^{1}-Y_{s}^{2},
f_{1}(s,Y_{s}^{1},Z_{s}^{1})-f_{2}(s,Y_{s}^{2},Z_{s}^{2})\rangle ds
\]
and
\[
I_2:=\dint_{t}^{\tau}\overline{e}^2_s\Delta_{s}^{\frac{\beta}{2}-2}
\sum_{j=1}^{r}
\left(\sum_{i=1}^{d}(Y_{i,s}^{1}-Y_{i,s}^{2})(Z_{i,j,s}^{1}-Z_{i,j,s}^{2})\right)^{2}ds.
\]
 In order to complete the proof of Lemma \ref{lem43} we need to
 estimate $I_1$ and $I_2$.
 \\\\
 \bf{Estimate of $I_1$.}  Let $\Phi(s):=|Y_{s}^{1}| + |Y_{s}^{2}|+
|Z_{s}^{1}| + |Z_{s}^{2}|+v_s.$ Since $\1_{\{\Phi_s \;\;\leq\;\; N\}} \leq
\1_{\{v_s \;\;\leq\;\; N\}}$ and $ f $ satisfies $\bf{(H4)}$, then a simple
computation shows that
\begin{align*}
\langle Y_{s}^{1}-Y_{s}^{2},&
f_{1}(s,Y_{s}^{1},Z_{s}^{1})-f_{2}(s,Y_{s}^{2},Z_{s}^{2})\rangle
\\&\vspace{0.4cm}\leq\;\overline{e}^{\frac{-1}{2}}_s\Delta^{\frac{1}{2}}_{s}
|f_{1}(s,Y_{s}^{1},Z_{s}^{1})-f_{2}(s,Y_{s}^{2},Z_{s}^{2})| \1_{\{\Phi_s
>N\}}\\ &\vspace{0.4cm}\quad  + 2N [\rho_N (f_1-f)_s+\rho_N (f_2-f)_s]\1_{\{v_s
\leq  N\}} \\ & \vspace{0.4cm}\quad+ [K"\log (A_N)\overline{e}^{-1}_s
\Delta_{s}+\overline{\gamma}\left| Z_{s}^{1}-Z_{s}^{2}\right|^{2}]\1_{\{\Phi_s
\;\;\leq\;\; N\}}
\end{align*}
\\ Therefore, using Lemma
\ref{lem40} with $h_s = \Delta_s$, we get
\begin{align*}
& I_1 \leq\;\dint_t^\tau\overline{e}^{\frac{1}{2}}_s\Delta^{\frac{\beta
-1}{2}}_{s} |f_{1}(s,Y_{s}^{1},Z_{s}^{1})-f_{2}(s,Y_{s}^{2},Z_{s}^{2})|
\1_{\{\Phi_s >N\}}ds  \\ & \;\;\; \;\;\;
+2N\dint_t^\tau\overline{e}_s\Delta^{\frac{\beta}{2}-1}_{s} [\rho_N
(f_1-f)_s+\rho_N (f_2-f)_s]\1_{\{v_s \leq N\}} ds \\ &\;\;\;
\;\;\;+\dint_t^\tau\overline{e}_s\Delta^{\frac{\beta}{2}-1}_{s} [K"\log
(A_N)\overline{e}^{-1}_s \Delta_{s}+\overline{\gamma}\left|
Z_{s}^{1}-Z_{s}^{2}\right|^{2}]\1_{\{\Phi_s \;\;\leq\;\; N\}}ds
\\
&\leq \varepsilon\sup_{s\in [u,u']}\Delta_{s}^{\frac{\beta}{2}}
\\&\;\;\;
\;\;\; + \varepsilon^{(1-\beta)}
\overline{e}^{\frac{\beta}{2}}_{u'} \dint_{u}^{u'}
|f_{1}(s,Y_{s}^{1},Z_{s}^{1})-f_{2}(s,Y_{s}^{2},Z_{s}^{2})|^\beta \1_{\{\Phi_s
>N\}} ds\\ &\;\;\; \;\;\; +
2N\dint_t^\tau\overline{e}_s\Delta^{\frac{\beta}{2}-1}_{s} [\rho_N
(f_1-f)_s+\rho_N (f_2-f)_s]\1_{\{v_s \leq N\}} ds \\ & \;\;\;
\;\;\;+\dint_t^\tau\overline{e}_s\Delta^{\frac{\beta}{2}-1}_{s} [K"\log
(A_N)\overline{e}^{-1}_s \Delta_{s}+\overline{\gamma}\left|
Z_{s}^{1}-Z_{s}^{2}\right|^{2}]\1_{\{\Phi_s \;\;\leq\;\; N\}}ds
\end{align*}
\bf{Estimate of $I_2$.} \ Since
$$\sum_{j=1}^{r}\left(\sum_{i=1}^{d}(Y_{i,s}^{1}-Y_{i,s}^{2})(Z_{i,j,s}^{1}-Z_{i,j,s}^{2})\right)^{2}
\leq \left| Y_{s}^{1}-Y_{s}^{2}\right|^{2}\left| Z_{s}^{1}-Z_{s}^{2}\right|^{2}
\leq \overline{e}^{-1}_s\Delta_{s}\left| Z_{s}^{1}-Z_{s}^{2}\right|^{2}$$ then
$$I_2\leq\dint_t^\tau\overline{e}_s\Delta^{\frac{\beta}{2}-1}_{s}\left|
Z_{s}^{1}-Z_{s}^{2}\right|^{2} ds. $$
%
%
\noindent Now, coming back to equation (\ref{equ1}) and taking into account the
above estimates we get for every $ \varepsilon
> 0$,
\begin{equation}\label{equ2}
 \begin{array}{lll}
&& \Delta_{t}^{\frac{\beta}{2}} +\dfrac{\beta}{2}
(\beta-1-2\overline{\gamma})\dint_{t}^{\tau}
\overline{e}_s\Delta_{s}^{\frac{\beta}{2}-1}
\left|Z_{s}^{1}-Z_{s}^{2}\right|^{2}ds
\\&&\\&&\quad\;
\;\;\leq\;\; \overline{e}^{\frac{\beta}{2}}_\tau |Y^1_\tau
-Y^2_\tau|^{\beta}+\dfrac{\overline{e}^
{\frac{\beta}{2}}_{u'}}{A^{\frac{\beta}{2}}_N}
 + \beta\varepsilon\displaystyle\sup_{s\in
[u,u']}\Delta_{s}^{\frac{\beta}{2}} \\ &&\\&&\quad\quad\,\,\,
+\beta\varepsilon^{(1-\beta)}
\overline{e}^{\frac{\beta}{2}}_{u'} \dint_{u}^{u'}
|f_{1}(s,Y_{s}^{1},Z_{s}^{1})-f_{2}(s,Y_{s}^{2},Z_{s}^{2})|^\beta \1_{\{\Phi_s
>N\}} ds
\\ &&\\ &&\quad\quad\,\,\,
+2N\beta\overline{e}^{\frac{\beta}{2}}_{\tau}
A_N^{1-\frac{\beta}{2}}\dint_{u}^{\tau}\rho_N (f_1-f)_s+\rho_N (f_2-f)_s
\1_{\{v_s \;\;\leq\;\; N\}} ds
\\ &&\\&&\quad\quad\,\,\,
-\beta\dint_{t}^{\tau}\overline{e}_s\Delta_{s}^{\frac{\beta}{2}-1}\langle
Y_{s}^{1}-Y_{s}^{2}, \left( Z_{s}^{1}-Z_{s}^{2}\right)dW_{s}\rangle.
\end{array}
\end{equation}
%
%
For a given $\hbar>1$, let $\tau_{\hbar}$ be the stopping time defined by $$
\tau_{\hbar}:= \inf\{s\geq u, \quad
|Y_s^1-Y_s^2|^2+\dint_u^s|Z^1_r-Z^2_r|^2dr\geq \hbar \}\wedge u',$$ Choose
$\tau=\tau_{\hbar}$, $ t=u $, \ then pass to the expectation in equation
(\ref{equ2}) to obtain, when $\hbar \rightarrow \infty$,
\begin{equation}\label{equ3}
\begin{array}{lll}
& &\dfrac{\beta}{2} (\beta - 1-2\overline{\gamma})\E\dint_{u}^{u'}
\overline{e}_s\Delta_{s}^{\frac{\beta}{2}-1}\left|Z_{s}^{1}-Z_{s}^{2}\right|^{2}ds
\\&&\\&&
\;\;\leq\;\; \overline{e}^{\frac{\beta}{2}}_{u'}
\E(|Y^1_{u'}-Y^2_{u'}|^{\beta}) + \dfrac{\overline{e}^{\frac{\beta}{2}}_{u'}}
{A^{\frac{\beta}{2}}_N} + \beta\varepsilon\E\displaystyle(\sup_{s\in
[u,u']}\Delta_{s}^{\frac{\beta}{2}}) \\&&\\&&\;\;\; \;\;\;
+\beta\varepsilon^{(1-\beta)}
\overline{e}^{\frac{\beta}{2}}_{u'} \E\dint_{u}^{u'}
|f_{1}(s,Y_{s}^{1},Z_{s}^{1})-f_{2}(s,Y_{s}^{2},Z_{s}^{2})|^\beta \1_{\{\Phi_s
>N\}} ds
\\ &&\\&&\;\;\;
\;\;\; +2N\beta\overline{e}^{\frac{\beta}{2}}_{u'}
A_N^{1-\frac{\beta}{2}}\E\dint_{u}^{u'}\rho_N (f_1-f)_s+\rho_N (f_2-f)_s
\1_{\{v_s \;\;\leq\;\; N\}} ds.
\end{array}
\end{equation}
%
%
Return back to (\ref{equ2}) and use the Burkholder-Davis-Gundy inequality to show that there exists a
universal constant $ c $ such that
\begin{align*} \E(\sup_{u \leq t \leq T} \Delta_{t}^{\frac{\beta}{2}})
&
\leq\; \overline{e}^{\frac{\beta}{2}}_{u'} \E (|Y^1_{u'}-Y^2_{u'}|^{\beta}) +
\dfrac{\overline{e}^{\frac{\beta}{2}}_{u'}}{A^{\frac{\beta}{2}}_N} +
\beta\varepsilon\E(\sup_{s\in [u,u']}\Delta_{s}^{\frac{\beta}{2}})
\\&\;\;\;
\;\;\; +\beta\varepsilon^{(1-\beta)}
\overline{e}^{\frac{\beta}{2}}_{u'} \E\dint_{u}^{u'}
|f_{1}(s,Y_{s}^{1},Z_{s}^{1})-f_{2}(s,Y_{s}^{2},Z_{s}^{2})|^\beta \1_{\{\Phi_s
>N\}} ds
\\ &\;\;\;
\;\;\; +2N\beta\overline{e}^{\frac{\beta}{2}}_{u'}
A_N^{1-\frac{\beta}{2}}\E\dint_{u}^{u'}\rho_N (f_1-f)_s+\rho_N (f_2-f)_s
\1_{\{v_s \;\;\leq\;\; N\}} ds
\\ &\;\;\;
\;\;\; +c\beta\E (
\dint_{u}^{T}\overline{e}^2_s\Delta_{s}^{\beta-2}\sum_{j=1}^r[\sum_{i=1}^d
(Y_{i,s}^{1}-Y_{i,s}^{2})(Z_{ij,s}^{1}-Z_{ij,s}^{2})]^2ds)^{\frac{1}{2}}.
\end{align*}
%
%
%
%
But, there exists a positive constant $ C_\beta $ depending only on $ \beta $
such that
\begin{align*}
 c\beta\E (
\dint_{u}^{u'}\overline{e}^2_s\Delta_{s}^{\beta-2} &\sum_{j=1}^r[\sum_{i=1}^d
(Y_{i,s}^{1}-Y_{i,s}^{2})(Z_{ij,s}^{1}-Z_{ij,s}^{2})]^2 ds)^{\frac{1}{2}} \\ &
\leq  \,\dfrac{1}{4}\E(\sup_{u\leq t\leq u'}
\Delta_{t}^{\frac{\beta}{2}}) + C_\beta\E
\dint_{u}^{u'}\overline{e}_s\Delta_{s}^{\frac{\beta}{2}-1}
|Z_{s}^{1}-Z_{s}^{2}|^2 ds. \end{align*}
%
%
Use (\ref{equ3}) and take $\varepsilon$ small enough to obtain the existence of
a positive constant $C=C(\beta,\overline{\gamma})$
  such that
\begin{align*}
& \E(\sup_{u\leq t\leq u'} \Delta_{t}^{\frac{\beta}{2}}) + \E
\dint_{u}^{u'}\overline{e}_s\Delta_{s}^{\frac{\beta}{2}-1}
|Z_{s}^{1}-Z_{s}^{2}|^2 ds
\\&
\;\;\leq\;\; C \left[\overline{e}^{\frac{\beta}{2}}_{u'} \E
|Y^1_{u'}-Y^2_{u'}|^{\beta}+\dfrac{\overline{e}^{\frac{\beta}{2}}_{u'}}
{A^{\frac{\beta}{2}}_N}\right.
+\overline{e}^{\frac{\beta}{2}}_{u'}\sup_{i}\E\dint_{u}^{u'}
|f_{i}(s,Y_{s}^{i},Z_{s}^{i})|^\beta \1_{\{\Phi_s >N\}} ds
\\ &\;\;\;
\;\;\; \left. +N\overline{e}^{\frac{\beta}{2}}_{u'}
A_N^{1-\frac{\beta}{2}}\E\dint_{u}^{u'}\rho_N (f_1-f)_s+\rho_N (f_2-f)_s
\1_{\{v_s \;\;\leq\;\; N\}} ds \right].
\end{align*}
We shall estimate $ J:=\sup_{i}\E\dint_{u}^{u'}
|f_{i}(s,Y_{s}^{i},Z_{s}^{i})|^\beta \1_{\{\Phi_s >N\}} ds$, $i=1,2$.
\\ Using
the fact that  $\1_{\{\Phi_s >N\}}\leq \1_{\{v_s
>5^{-1}N\}}+\1_{\{|Y^1_s| >5^{-1}N\}}+\1_{\{|Y^2_s| >5^{-1}N\}} +
\1_{\{|Z^1_s| >5^{-1}N\}}+\1_{\{|Z^2_s| >5^{-1}N\}} $ and  $ \1_{\{a > b\}}\leq
\dfrac{a^\nu}{b^\nu}$ for every $ a, b,\nu >0 $, we show that for  every  $
N>1$
\begin{align*}
J\leq & \; \left(\dfrac{5}{N}\right)^{\nu}
 \sup_{i}\E\dint_{u}^{u'}
|f_{i}(s,Y_{s}^{i},Z_{s}^{i})|^\beta v^{\nu}_s ds  \\ & +
  \left(\dfrac{5}{N}\right)^{\nu}
 \sup_{i}\E\dint_{u}^{u'}
|f_{i}(s,Y_{s}^{i},Z_{s}^{i})|^\beta |Y^1_s|^{\nu} ds \\ &  +
  \left(\dfrac{5}{N}\right)^{\nu}
 \sup_{i}\E\dint_{u}^{u'}
|f_{i}(s,Y_{s}^{i},Z_{s}^{i})|^\beta |Y^2_s|^{\nu} ds \\ & +
\left(\dfrac{5}{N}\right)^{\nu}
 \sup_{i}\E\dint_{u}^{u'}
|f_{i}(s,Y_{s}^{i},Z_{s}^{i})|^\beta |Z^1_s|^{\nu} ds. \\ & +
\left(\dfrac{5}{N}\right)^{\nu}
 \sup_{i}\E\dint_{u}^{u'}
|f_{i}(s,Y_{s}^{i},Z_{s}^{i})|^\beta |Z^2_s|^{\nu}
 ds.
\end{align*}
%
%
 using Young's inequality, one can prove that there exists
 a positive constant $C $ such that for every
 $ N>1$
\begin{align*}
& J \leq   \dfrac{C}{N^{\nu}} \left\{ 1+ \Theta^1_p + \Theta^2_p
+\sup_{i}\E\dint_{u}^{u'} |f_{i}(s,Y_{s}^{i},Z_{s}^{i})|^{\beta
(\frac{q'}{q'-\nu}\vee\frac{2}{2-\nu}\vee\frac{p}{p-\nu})}
 ds + \E\dint_{u}^{u'}
v_s^{q'} ds\right\}.
\end{align*}
\\
where $\Theta^i_p:= \E(\sup_t|Y^i_t|^p)+\E\left(\dint_{0}^{T}|Z^i_{s}|^{2} ds
\right)^{\frac{p}{2}}$. \\
 Using Proposition \ref{pro42},  we get
(since $ \beta (\frac{q'}{q'-\nu}\vee\frac{2}{2-\nu}\vee\frac{p}{p-\nu})\leq
\hat{\beta}$)
\begin{align*}
& J \leq   \dfrac{C}{N^{\nu}} \left\{ 1+ \Theta^1_p + \Theta^2_p
+\E\dint_{0}^{T} |\overline{\eta}_s|^{q}
 ds + \E\dint_{u}^{u'}
v_s^{q'} ds\right\}.
\end{align*}
Hence, for $ a :=(\frac{\nu}{\mu}\wedge\frac{\beta}{2}) - \beta rK"$ and $N$
large enough we get (since $ A_N \leq N^\mu$ by assumption {\{bf{(H.4)}}),
\begin{align*}
& \E\sup_{u\leq t\leq u'} \Delta_{t}^{\frac{\beta}{2}}+\E
\dint_{u}^{u'}\overline{e}_s\Delta_{s}^{\frac{\beta}{2}-1}
|Z_{s}^{1}-Z_{s}^{2}|^2 ds
\\&
\;\;\leq\;\; N A_N^{1+\frac{\beta}{2}}\left[ \E
|Y^1_{u'}-Y^2_{u'}|^{\beta}+\E\dint_{0}^{T}\rho_N (f_1-f)_s+\rho_N (f_2-f)_s
\1_{\{v_s \;\;\leq\;\; N\}} ds \right]\\ &\;\;\; +
\dfrac{1}{A^{a}_N}\left[1+\Theta_p^1+\Theta_p^2+\E\dint_0^T
\overline{\eta}_s^qds+ \E\dint_0^T v_s^{q'}ds\right].
\end{align*}
 Lemma
\ref{lem43} is proved. \eop

As a consequence of lemma \ref{lem43}, we have


\begin{lemma}\label{lem44}
Let $(\xi^i, f_i)_{i=1,2}$ satisfies \bf{(H.3)} (with the same
$\overline{\eta}, \alpha$ and $\alpha'$) and let $(Y^i, Z^i)$ be a solution of
$(E^{(\xi^i, f_i)})$. Then, there exists $\beta=\beta(p,q,\alpha,\alpha') \in
]1,p\wedge 2[$
 such that for every $\varepsilon
> 0$
there is an integer
$N_\varepsilon=N_\varepsilon(p,q,\alpha,\alpha',K',\mu,q',\varepsilon,
(A_N)_N)$ such that for every function $f$ satisfying $\bf{(H.4)}$
\begin{align*}
& \E(\sup_{0\leq t\leq T} |Y^1_t-Y^2_t|^{\beta})+\E
\dint_{0}^{T}\dfrac{|Z_{s}^{1}-Z_{s}^{2}|^2}
{\left(1+|Y^1_s-Y^2_s|^2\right)^{1-\frac{\beta}{2}}}
 ds
\\&
 \leq  \;\; N_\varepsilon\left[ \E
|\xi^1-\xi^2|^{\beta}+\E\dint_{0}^{T}\rho_{N_\varepsilon}
(f_1-f)_s+\rho_{N_\varepsilon} (f_2-f)_s  ds \right] \\ &\;\;\;\;+ \varepsilon
 \left[1+\Theta_p^1+\Theta_p^2+\E\dint_0^T
\overline{\eta}_s^qds+ \E\dint_0^T v_s^{q'}ds\right].
\end{align*}
\end{lemma}


\noindent\bf{Proof.} Let $(u_0=0 <...< u_{\ell+1}=T)$ be a subdivision of $ [0,T]$ such that
for every  $ i\in\{0,..,\ell\} $ $$u_{i+1}-u_{i}\leq r$$ \noindent From lemma
\ref{lem43} we have : for all $\varepsilon
> 0$ there is an integer $N_\varepsilon$ such that for every
function $ f $ satisfying $\bf{(H.4)}$
\begin{align*}
&\E(\sup_{u_\ell\leq t\leq T} |Y^1_t-Y^2_t|^{\beta})+\E
\dint_{u_\ell}^{T}\dfrac{|Z_{s}^{1}-Z_{s}^{2}|^2}
{\left(1+|Y^1_s-Y^2_s|^2\right)^{1-\frac{\beta}{2}}}
 ds
\\&
 \leq  \;\; N_\varepsilon\left[ \E
(|\xi^1-\xi^2|^{\beta})+\E\dint_{0}^{T}\rho_{N_\varepsilon}
(f_1-f)_s+\rho_{N_\varepsilon} (f_2-f)_s  ds \right] \\ &\;\;\;\;+ \varepsilon
 \left[1+\Theta_p^1+\Theta_p^2+\E\dint_0^T
\overline{\eta}_s^qds+ \E\dint_0^T v_s^{q'}ds\right].
\end{align*}
Assume that for some $i\in \{0,..,\ell\}$ we have for all $\varepsilon
> 0$ there is an integer $N_\varepsilon$ such that for every
function  $ f $ satisfying $\bf{(H.4)}$
\begin{align*}
& \E(\sup_{u_{i+1}\leq t\leq T} |Y^1_t-Y^2_t|^{\beta})+\E
\dint_{u_{i+1}}^{T}\dfrac{|Z_{s}^{1}-Z_{s}^{2}|^2}
{\left(1+|Y^1_s-Y^2_s|^2\right)^{1-\frac{\beta}{2}}}ds
\\&
 \leq  \;\; N_\varepsilon\left[ \E
(|\xi^1-\xi^2|^{\beta})+\E\dint_{0}^{T}\rho_{N_\varepsilon}
(f_1-f)_s+\rho_{N_\varepsilon} (f_2-f)_s  ds \right]  \\ &\;\;\;\;+ \varepsilon
 \left[1+\Theta_p^1+\Theta_p^2+\E\dint_0^T
\overline{\eta}_s^qds+ \E\dint_0^T v_s^{q'}ds\right].
\end{align*}
Then, for every $\varepsilon'
> 0$ there is an integer $N_{\varepsilon'}$ such that for every
function  $ f $ satisfying $\bf{(H.4)}$
\begin{align*}
& \E(\sup_{u_{i}\leq t\leq T} |Y^1_t-Y^2_t|^{\beta})+\E
\dint_{u_{i}}^{T}\dfrac{|Z_{s}^{1}-Z_{s}^{2}|^2}
{\left(1+|Y^1_s-Y^2_s|^2\right)^{1-\frac{\beta}{2}}}ds
\\ & \leq \E(\sup_{u_{i}\leq t\leq u_{i+1}} |Y^1_t-Y^2_t|^{\beta})+\E
\dint_{u_{i}}^{u_{i+1}}\dfrac{|Z_{s}^{1}-Z_{s}^{2}|^2}
{\left(1+|Y^1_s-Y^2_s|^2\right)^{1-\frac{\beta}{2}}}ds
\\&
  \;\;+ N_{\varepsilon'}\left[ \E
(|\xi^1-\xi^2|^{\beta})+\E\dint_{0}^{T}\rho_{N_{\varepsilon'}}
(f_1-f)_s+\rho_{N_{\varepsilon'}} (f_2-f)_s  ds \right] \\ &\;\;\;\;+
\varepsilon'
 \left[1+\Theta_p^1+\Theta_p^2+\E\dint_0^T
\overline{\eta}_s^qds+ \E\dint_0^T v_s^{q'}ds\right].
\end{align*}
Using Lemma \ref{lem43} we obtain; for every $\varepsilon',\varepsilon"
> 0$ there exist $N_{\varepsilon'}>0$ and $N_{\varepsilon"}>0$ such that for every
function  $ f $ satisfying $\bf{(H.4)}$
\begin{align*}
& \E(\sup_{u_{i}\leq t\leq T} |Y^1_t-Y^2_t|^{\beta})+\E
\dint_{u_{i}}^{T}\dfrac{|Z_{s}^{1}-Z_{s}^{2}|^2}
{\left(1+|Y^1_s-Y^2_s|^2\right)^{1-\frac{\beta}{2}}}ds
\\ & \leq N_{\varepsilon"}\left[ \E
(|Y^1_{u_{i+1}}-Y^2_{u_{i+1}}|^{\beta})+\E\dint_{0}^{T}
\rho_{N_{\varepsilon"}}
(f_1-f)_s+\rho_{N_{\varepsilon"}} (f_2-f)_s  ds \right] \\ &
  \;\;+ N_{\varepsilon'}\left[ \E
(|\xi^1-\xi^2|^{\beta})+\E\dint_{0}^{T}\rho_{N_{\varepsilon'}}
(f_1-f)_s+\rho_{N_{\varepsilon'}} (f_2-f)_s  ds \right] \\ &\;\;\;\;+
2\varepsilon'
 \left[1+\Theta_p^1+\Theta_p^2+\E\dint_0^T
\overline{\eta}^qds+ \E\dint_0^T v_s^{q'}ds\right]
 \\ & \leq  N_{\varepsilon'}N_{\varepsilon"}\;\;\E
(|\xi^1-\xi^2|^{\beta}) \\ & \;\; +
(N_{\varepsilon'}N_{\varepsilon"}+2N_{\varepsilon'})\;\;\E\dint_{0}^{T}\rho_{(N_{\varepsilon'}N_{\varepsilon"})}
(f_1-f)_s+\rho_{(N_{\varepsilon'}N_{\varepsilon"})} (f_2-f)_s  ds
\\ &\;\;\;\;+ (2\varepsilon'+\varepsilon"N_{\varepsilon'})
 \left[1+\Theta_p^1+\Theta_p^2+\E\dint_0^T
\overline{\eta}_s^qds+ \E\dint_0^T v_s^{q'}ds\right].
\end{align*} For $\varepsilon > 0$, let
$\varepsilon':=\dfrac{\varepsilon}{4}$ and $\varepsilon":=
\dfrac{\varepsilon}{2N_{(\frac{\varepsilon}{4})}}$, we then deduce that there  exists
 an integer $N_\varepsilon$ such that for every function  $ f $ satisfying
$\bf{(H.4)}$
\begin{align*}
& \E(\sup_{u_{i}\leq t\leq T} |Y^1_t-Y^2_t|^{\beta})+\E
\dint_{u_{i}}^{T}\dfrac{|Z_{s}^{1}-Z_{s}^{2}|^2}
{\left(1+|Y^1_s-Y^2_s|^2\right)^{1-\frac{\beta}{2}}}ds
\\&
 \leq  \;\; N_\varepsilon\left[ \E
(|\xi^1-\xi^2|^{\beta})+\E\dint_{0}^{T}\rho_{N_\varepsilon}
(f_1-f)_s+\rho_{N_\varepsilon} (f_2-f)_s  ds \right]\\ &\;\;\;\;+ \varepsilon
 \left[1+\Theta_p^1+\Theta_p^2+\E\dint_0^T
\overline{\eta}_s^qds+ \E\dint_0^T v_s^{q'}ds\right].
\end{align*}
We complete the proof by induction \eop \vskip1cm


\begin{proposition}\label{pro43}
Let $(\xi^i, f_i)_{i=1,2}$ satisfies \bf{(H.3)} (with the same
$\overline{\eta}, \alpha$ and $\alpha'$) and let $(Y^i, Z^i)$ be a solution of
$(E^{(\xi^i, f_i)})$. Then, there exists $\beta=\beta(p,q,\alpha,\alpha') \in
]1,p\wedge 2[$
 such that for every $\varepsilon
> 0$
there is an integer
$N_\varepsilon=N_\varepsilon(p,q,\alpha,\alpha',K',\mu,q',\varepsilon,
(A_N)_N)$ such that for every function $f$ satisfying $\bf{(H.4)}$
\begin{align*}
& \E(\sup_{0\leq t\leq T} |Y^1_t-Y^2_t|^{\beta})+\E
\left(\dint_{0}^{T}|Z_{s}^{1}-Z_{s}^{2}|^2ds\right)^\frac{\beta}{2}
\\&
 \leq  \;\; N_\varepsilon\left[ \E
(|\xi^1-\xi^2|^{\beta})+\E\dint_{0}^{T}\rho_{N_\varepsilon}
(f_1-f)_s+\rho_{N_\varepsilon} (f_2-f)_s  ds \right] \\ &\;\;\;\;+ \varepsilon
 \left[1+\Theta_p^1+\Theta_p^2+\E\dint_0^T
\overline{\eta}_s^qds+ \E\dint_0^T v_s^{q'}ds\right],
\end{align*}
where \ $ \Theta^i_p:= \E(\sup_t|Y^i_t|^p)+\E\left(\dint_{0}^{T}|Z^i_{s}|^{2} ds
\right)^{\frac{p}{2}}.$
\end{proposition}

\noindent\bf{Proof.} Using H\"{o}lder's inequality, Young's inequality and the fact that
$\dfrac{\beta}{2}<1$, we obtain for all $\varepsilon'
>0$

\begin{align*}
& \E \left(\dint_{0}^{T}|Z_{s}^{1}-Z_{s}^{2}|^2ds\right)^\frac{\beta}{2}
\\&
 \leq  \;\; \E\;\big\{
\left[\dint_{0}^{T}\dfrac{|Z_{s}^{1}-Z_{s}^{2}|^2}
{\left(1+|Y^1_s-Y^2_s|^2\right)^{1-\frac{\beta}{2}}}
 ds\right]^\frac{\beta}{2}
\sup_{s\leq T}
\left(1+|Y^1_s-Y^2_s|^2\right)^{(1-\frac{\beta}{2})\frac{\beta}{2}}\big\}
\\&
 \leq  \;\;
\left[\E\dint_{0}^{T}\dfrac{|Z_{s}^{1}-Z_{s}^{2}|^2}
{\left(1+|Y^1_s-Y^2_s|^2\right)^{1-\frac{\beta}{2}}}
 ds\right]^\frac{\beta}{2}
\left(1+\E(\sup_{s\leq T} |Y^1_s-Y^2_s|^\beta)\right)^{\frac{2-\beta}{2}}
\\&
 \leq  \;\;\left[\E(\sup_{s\leq T}
|Y^1_s-Y^2_s|^\beta) + \E\dint_{0}^{T}\dfrac{|Z_{s}^{1}-Z_{s}^{2}|^2}
{\left(1+|Y^1_s-Y^2_s|^2\right)^{1-\frac{\beta}{2}}}
 ds \right]^\frac{\beta}{2}
 \\ & \;\;+ \left[\E(\sup_{0\leq t\leq T} |Y^1_t-Y^2_t|^{\beta})+\E
\dint_{0}^{T}\dfrac{|Z_{s}^{1}-Z_{s}^{2}|^2}
{\left(1+|Y^1_s-Y^2_s|^2\right)^{1-\frac{\beta}{2}}}
 ds\right] \\&
 \leq  \;\;\varepsilon' \; + \;(1+\varepsilon^{'\frac{\beta - 2}{\beta}})
\left[\E(\sup_{0\leq t\leq T} |Y^1_t-Y^2_t|^{\beta})+\E
\dint_{0}^{T}\dfrac{|Z_{s}^{1}-Z_{s}^{2}|^2}
{\left(1+|Y^1_s-Y^2_s|^2\right)^{1-\frac{\beta}{2}}}
 ds\right] .
\end{align*}

Use lemma \ref{lem43} to conclude that for every $\varepsilon', \varepsilon">0$
\begin{align*}
& \E \left(\dint_{0}^{T}|Z_{s}^{1}-Z_{s}^{2}|^2ds\right)^\frac{\beta}{2}
\\&
  \leq  \;\;\varepsilon' \; + \;(1+\varepsilon^{'\frac{\beta - 2}{\beta}})
  N_{\varepsilon"}\left[ \E
(|\xi^1-\xi^2|^{\beta})+\E\dint_{0}^{T}\rho_{N_{\varepsilon"}}
(f_1-f)_s+\rho_{N_{\varepsilon"}} (f_2-f)_s  ds \right] \\ &\;\;\;\;+
\varepsilon" (1+\varepsilon^{'\frac{\beta - 2}{\beta}})
 \left[1+\Theta_p^1+\Theta_p^2+\E\dint_0^T
\overline{\eta}_s^qds+ \E\dint_0^T v_s^{q'}ds\right].
\end{align*}
Letting $\varepsilon' = \dfrac{\varepsilon}{2}$ and $\varepsilon"
=
\dfrac{\varepsilon}{2(1+(\frac{\varepsilon}{2})^{\frac{\beta-2}{2}})}$, we
finish this proof of proposition \ref{pro43}. \eop


\vskip0.3cm
\begin{remark} The uniqueness of equation
$(E^{(\xi, f)})$ follows by letting $f_1 = f_2 = f$ and $\xi_1 = \xi_2 = \xi$
in Proposition \ref{pro43}.
\end{remark}

The following stability result follows from propositions (\ref{pro43}),
(\ref{pro42}) and (\ref{pro41})


\begin{proposition} \label{pro44} Let $ (\xi, f)$ satisfies
\bf{(H.0)-(H.4)} and $ (\xi^n, f_n)_n$  satisfies \bf{(H.0)-(H.3)} uniformly on
$n$. Assume moreover that

\vskip 0.15cm\noindent (a) $ \xi^n\rightarrow\xi \;\; \hbox{a.s.}$\, \ \  and \ \ $\sup_n\E\big(\vert\xi_n\vert^p\exp(\frac{p}{2}\int_0^T\lambda_sds)\big) < \infty$

\vskip 0.15cm\noindent (b) For every $N\in\N^*$, \ $\lim_n\rho_N(f_n-f)=0 \;\;\hbox{a.e.}$

\vskip 0.15cm\noindent (c) for every $n\in\N^*$, the BSDE \  $ (E^{(\xi^n, f_n)})$ \ has a solution $(Y^n,Z^n)$ which
satisfies,

\vskip 0.15cm \ \ \ \ \ \ \ \ \ $ \sup_n\E(\sup_{t\leq T}|Y^n_t|^pe^{\frac{p}{2}\int_0^T
\lambda_sds})<\infty$.

\vskip 0.15cm\noindent Then, there exists $ (Y, Z)\in\;\;\L^{p}(\Omega; \mathcal{C}([0,T];
\R^d))\times\L^{p}(\Omega; \L^{2}([0,T]; \R^{dr}))$  such that

\ni $i) $ \ \
$\displaystyle
  \E(\sup_t\mid Y_t\mid^p  e^{\frac{p}{2}\int_0^t \lambda_sds}) + \E\left[\dint_0^T
e^{\int_0^s \lambda_rdr}\mid Z_s\mid^2ds\right]^\frac{p}{2} $

\hskip 2.5cm $\displaystyle \leq
C^{p,\gamma} \left\{\E(\mid \xi\mid^pe^{\frac{p}{2}\int_0^T \lambda_sds})
+\E\left(\dint_0^T e^{\int_0^s \lambda_rdr}\eta_s ds\right)^{\frac{p}{2}}+
\E\left(\dint_0^T e^{\frac{1}{2}\int_0^s \lambda_rdr}f^0_s ds\right)^{p}
\right\}
$
\\ $ii)$\;\;  for every $p'<p$,  $(Y^n,Z^n)\longrightarrow (Y,
Z)\;\;\hbox{strongly in}\;\;\L^{p'}(\Omega; \mathcal{C}([0,T];
\R^d))\times\L^{p'}(\Omega; \L^{2}([0,T]; \R^{dr})). $

\vskip 0.15cm \noindent  $iii)$\;\; for every
$\hat{\beta}<\dfrac{2}{\alpha'}\wedge\dfrac{p}{\alpha}\wedge\dfrac{p}{\alpha'}\wedge
q$, \ \
$\displaystyle \lim_{n\rightarrow\infty}\E\dint_0^T
|f_n(s,Y^n_s,Z^n_s)-f(s,Y_s,Z_s)|^{\hat{\beta}}ds=0$

\vskip 0.15cm \noindent Moreover, $(Y,Z)$ is the unique solution of $(E^{(\xi,f)})$.
\end{proposition}


\noindent\bf{Proof.} From Proposition \ref{pro41}, Proposition \ref{pro42} and Proposition
\ref{pro43}, we have
 \vskip 0.1cm\ni  $ \begin{array}{l}
 a') \ \ \left[ \E(\sup_t|Y^n_t|^pe^{\frac{p}{2}\int_0^t
 \lambda_sds})
 +\E\left(\dint_{0}^{T}e^{\int_0^t \lambda_sds}|Z_{s}^{n}|^{2}
ds \right)^{\frac{p}{2}}\right]
\\ \ \ \ \ \ \ \ \ \leq \ C^{p,\gamma} \sup_n\left\{\E(\mid
\xi^n\mid^pe^{\frac{p}{2}\int_0^T \lambda_sds}) +\E\left(\dint_0^T e^{\int_0^s
\lambda_rdr}\eta_s ds\right)^{\frac{p}{2}}+ \E\left(\dint_0^T
e^{\frac{1}{2}\int_0^s \lambda_rdr}f^0_s ds\right)^{p} \right\} \\
\ \ \ \ \ \ \ := \ D.
 \end{array}$

\vskip 0.25cm\ni  $\; b') \ \ \E\dint_0^T |f_n(s, Y^n_s, Z^n_s)|^{\hat{\beta}}ds \ \leq \ C
(1+D+\dint \bar{\eta}^q_sds). $

\vskip 0.25cm\ni $\; c' ) \ \ $  There exists $\beta>1$ such that for every $\varepsilon>0$
there exists $N_\varepsilon>0$:
  \begin{align*} \E(\sup_t|Y^n_t-
 Y^m_t|^\beta)
 +\E\left(\dint_{0}^{T}|Z_{s}^{n}-Z_{s}^{m}|^{2}
ds \right)^{\frac{\beta}{2}} & \leq N_\varepsilon\E\left[
|\xi^n-\xi^m|^{\beta}+\dint_{0}^{T}\rho_{N_\varepsilon}
(f_n-f)_s+\rho_{N_\varepsilon} (f_m-f)_s  ds \right] \\  &\quad + \varepsilon
 \left[1+2D +\E\dint_0^T
\overline{\eta}_s^qds+ \E\dint_0^T v_s^{q'}ds\right].
 \end{align*}
\ni We deduce the existence of $ (Y, Z)\in\;\;\L^{p}(\Omega; \mathcal{C}([0,T];
\R^d))\times\L^{p}(\Omega; \L^{2}([0,T]; \R^{dr}))$  such that

\ni i) \ \ $\displaystyle  \E(\sup_t\mid Y_t\mid^p e^{\frac{p}{2}\int_0^t
\lambda_sds}) + \E\left[\dint_0^T e^{\int_0^s \lambda_rdr}\mid
Z_s\mid^2ds\right]^\frac{p}{2}$

 \hskip 2cm $\displaystyle \leq  C^{p,\gamma} \left\{\E(\mid
\xi\mid^pe^{\frac{p}{2}\int_0^T \lambda_sds}) +\E\left(\dint_0^T e^{\int_0^s
\lambda_rdr}\eta_s ds\right)^{\frac{p}{2}}+ \E\left(\dint_0^T
e^{\frac{1}{2}\int_0^s \lambda_rdr}f^0_s ds\right)^{p} \right\}
$
\\ $ii)$\;\;  for all $p'<p$,  $(Y^n,Z^n)\longrightarrow (Y,
Z)\;\;\hbox{strongly in}\;\;\L^{p'}(\Omega; \mathcal{C}([0,T];
\R^d))\times\L^{p'}(\Omega; \L^{2}([0,T]; \R^{dr})). $

\vskip 0.2cm\noindent
  Let us prove $iii)$. \ Set $a :=
\limsup_{n\rightarrow\infty}\E\dint_0^T
|f(s,Y^n_s,Z^n_s)-f(s,Y_s,Z_s)|^{\hat{\beta}}ds$. Consider a subsequence
$n'$ of $n$ such that $a := \lim_{n'\rightarrow\infty}\E\dint_0^T
|f(s,Y^{n'}_s,Z^{n'}_s)-f(s,Y_s,Z_s)|^{\hat{\beta}}ds$ and,
 $(Y^{n'},Z^{n'})\rightarrow (Y, Z)$ a.e. \\
 Assumption \bf{(H.3)} and the continuity of $f$
  ensure that $a=0$. It remains to prove that

$$ \limsup_{n\rightarrow\infty}\E\dint_0^T
|f_n(s,Y^n_s,Z^n_s)-f(s,Y^n_s,Z^n_s)|^{\hat{\beta}}ds=0$$
We use H$\ddot{o}$lder's
inequality, the previous claim b'), Proposition \ref{pro42} and Chebychev's
inequality to get $$
\begin{array}{l}\E\dint_0^T |f_n(s,Y^n_s,Z^n_s)-f(s,Y^n_s,Z^n_s)|^{\hat{\beta}}ds\\
\leq \E\dint_0^T \rho_N(f_n-f)_s^{\hat{\beta}}ds +  (\E\dint_0^T
|f_n(s,Y^n_s,Z^n_s)-f(s,Y^n_s,Z^n_s)|^{r\hat{\beta}}ds)^{\frac{1}{r}}
(\E\dint_0^T 1_{|Y^n_s|+|Z^n_s|\geq N}ds)^\frac{r-1}{r} \\ \leq\E\dint_0^T
\rho_N(f_n-f)_s^{\hat{\beta}}ds +
\dfrac{C(r)}{N^{\frac{(r-1)(p\wedge2)}{r}}},
\end{array}$$
for some reel $r>1$ such that
$r\hat{\beta}<\dfrac{2}{\alpha'}\wedge\dfrac{p}{\alpha}\wedge\dfrac{p}{\alpha'}\wedge
q.$

\ni We successively let $n\longrightarrow\infty$ and  $N\longrightarrow\infty$ to derive
 assertion $iii).$
 Proposition \ref{pro44} is proved \eop

 %
 %

\subsection{\bf{Approximation}} We shall construct a
sequence  $(\xi^n, f_n)$ which converges in a suitable sense to $(\xi,f)$ and which has good properties. With the help of this approximation, we can construct a
solution $(Y, Z)$ to the BSDE $(E^{(\xi,f)})$ by using Proposition  \ref{pro44}. \\
\\ \ni Let $h_t$ is a
predictable process such that $ 0<h_t\leq 1$ and set $ \overline{\Lambda}_t := \eta_t
+\overline{\eta}_t+f^0_t+M_t+K_t +\frac{1}{h_t}$ \ \
\begin{proposition}
\label{pro45} {\it Assume that $(\xi,f)$ satisfies} \bf{(H.0)}--\bf{(H.3)}.
{\it Then there exists a sequence
 $(\xi^n,f_n)$ such that
\par\noindent
$(a)$ \ For each $n$, $\xi^n$ is bounded, $|\xi^n|\leq |\xi|$ and $\xi^n$
converges to $\xi$ $a.s$.
\\ $(b)$ \; For each $n$, $f_n$ is uniformly Lipschitz in $(y,z)$.
\\$(c)$ \; $\vert f_n(t,\omega,y,z)\vert\; \;\;\leq\;\; \1_{\{
\overline{\Lambda}_t\;\;\leq\;\; n,\,\, \mid y\mid\leq n, \,\,
\mid z\mid\leq n\}} \{ \overline{\eta}_t + \mid y\mid^{\alpha}+
\mid z\mid^{\alpha'}+ 2p h_t\}\leq 2p +3n^{p}$.
\\
$(d)$ \; $  <y,f_n(t,\omega,y,z)> \;\;\leq\;\; \1_{\{
\overline{\Lambda}_t\;\;\leq\;\; n\}} \{\eta_t+ f^0_t \vert y\vert+ M_t \vert
y\vert^2 + K_t\vert y\vert\vert z\vert +
 10 h_t\}.
 $
 \\
 $(e)$ \; For every $N$, \ $\rho_N
(f_n-f)(t,\omega)\longrightarrow 0$ as $n\longrightarrow\infty$ a.e
$(t,\omega)$.
\\
$(f)$ \; For every $N$, \ $\rho_N (f_n-f)(t,\omega)\;\;\leq\;\; 2
\{\overline{\eta}_t + N^{\alpha}+  N^{\alpha'}+ 2p h_t\} $.}
\end{proposition}


\noindent\bf{Proof.} Let $ \psi: \R\longrightarrow [0,\dfrac{\exp(-1)}{c_1}] $ defined by: $$
\psi(x):= \left\{
\begin{tabular}{ll}
$ c_1^{-1}\exp{(-\dfrac{1}{1-x^2})} $ & if $|x|<1$  \\ 0 & else
\end{tabular} \right.
$$ where $ c_1= \dint_{-1}^1\exp{(-\dfrac{1}{1-x^2})} dx $.\\

\noindent Let  $ m := \dfrac{n^{2p}}{ h_t} $.\ the sequence $ (\xi^n,f_n)$ defined by : \ $\xi^n := \xi 1_{\{\mid
\xi\mid\leq n\}} $ and
\begin{align*}
 f_{n}(t,y,z) = &(c_1e)^2 \1_{\{\overline{\Lambda}_t\;\;\leq\;\;
n\}}\psi(n^{-2}|y|^2)\psi(n^{-2}|z|^2) \times \\ & \ \
 m^{(d+dr)}\dint_{\R^d}\int_{\R^{dr}} f(t,y-u,z-v)\Pi_{i=1}^d\psi
(m u_i)\Pi_{i=1}^d \Pi_{j=1}^r\psi(m v_{ij})dudv,
\end{align*}
satisfies the required properties. Indeed, (a) \ is obvious. $(e)$ \ follows from the definition of $f_n$.
$(f)$ \ follows from assumption $\bf{(H.3)}$ and assertion $(c)$. We shall prove assertions (b), (c) and (d).

$(b)$ \ For a fixed $t$ and $\omega$, $f_n(t,\omega, .,.)$ is smooth and
 with compact support in $[-n, n]^{d+dr}$. Moreover
$$
\mid \nabla_{y, z} f_n(t, \omega, y,z) \mid \leq C n^{2p +2},
$$
where $\nabla$ denotes the gradient and $C$ is a positive constant.

$(c)$ \ For all $(t, \omega, y, z)$ such that
$\overline{\Lambda}_t\;\;\leq\;\; n$, \,\, $\mid y\mid\leq n $ and
$\mid z\mid\leq n$ we obtain, by using assumption $\bf{(H.3)}$, that
\begin{align*}
 \mid f_{n}(t,y,z)\mid  \leq  &
 m^{(d+dr)}\dint_{\R^d}\int_{\R^{dr}} \mid f(t,y-u,z-v)\mid\Pi_{i=1}^d\psi
(m u_i)\Pi_{i=1}^d \Pi_{j=1}^r\psi(m v_{ij})dudv
\\ & \leq \overline{\eta}_t + \mid y\mid^{\alpha}+
\mid z\mid^{\alpha'} + m^{d}\dint_{\R^d} \bigg(\mid
y-u\mid^{\alpha} -\mid y\mid^{\alpha}\bigg)\Pi_{i=1}^d\psi (m
u_i)du\\ &  + m^{dr}\int_{\R^{dr}} \bigg(\mid z-v\mid^{\alpha'}
-\mid z \mid^{\alpha'}\bigg)\Pi_{i=1}^d \Pi_{j=1}^r\psi(m
v_{ij})dv
\\ & \leq  \overline{\eta}_t + \mid y\mid^{\alpha}+
\mid z\mid^{\alpha'}+\alpha (n +\frac{h_t}{n^{2p}})^{\alpha
-1}\frac{h_t}{n^{2p}} +\alpha' (n +\frac{h_t}{n^{2p}})^{\alpha'
-1}\frac{h_t}{n^{2p}}
\\ & \leq  \overline{\eta}_t + \mid y\mid^{\alpha}+
\mid z\mid^{\alpha'}+ 2p h_t
\end{align*}

$(d)$ \ For all $(t, \omega, y, z)$ such that
$\overline{\Lambda}_t\;\;\leq\;\; n$, \,\, $\mid y\mid\leq n $ and
$\mid z\mid\leq n$ \ we obtain, by using assumptions
$\bf{(H.2)}-\bf{(H.3)}$, that
\begin{align*}
 \langle y, f_{n}(t,y,z)\rangle  \leq  &
(c_1e)^2\psi(n^{-2}|y|^2)\psi(n^{-2}|z|^2)\times \\&
m^{(d+dr)}\dint_{\R^d}\int_{\R^{dr}} \langle f(t,y-u,z-v), y
-u\rangle\Pi_{i=1}^d\psi (m u_i)\Pi_{i=1}^d \Pi_{j=1}^r\psi(m
v_{ij})dudv \\ & + m^{(d+dr)}\dint_{\R^d}\int_{\R^{dr}} \mid
f(t,y-u,z-v)\mid \mid u\mid\Pi_{i=1}^d\psi (m u_i)\Pi_{i=1}^d
\Pi_{j=1}^r\psi(m v_{ij})dudv
\\ & \leq \eta_t+ f^0_t \vert y\vert+ M_t \vert
y\vert^2 + K_t\vert y\vert\vert z\vert +
 10 h_t
\end{align*}
 \eop
\begin{remark}
Theorem \ref{the21} follows now from Proposition \ref{pro44} and Proposition \ref{pro45}.
\end{remark}

%
%

\section{\bf{Application to partial differential equations (PDEs)}}
In this section, we consider the system of semilinear PDEs associated to the Markovian version of the BSDE $(E^{\xi, f)})$, for which we establish the existence and uniqueness of a weak (Sobolev) solution. In particular, we give a new feature which consists to prove, by using BSDEs techniques,
that the uniqueness for a nonhomogeneous system of semilinear PDE follows from the uniqueness of its associated homogeneous system of linear PDE.

\subsection{\bf{Formulation of the problem.}}
Let $\sigma : \R^k\longmapsto\R^{kr}$, \ $b : \R^k\longmapsto\R^{k}$ , \ $g : \R^k\longmapsto\R^{k}$, \ and  $F :
[0, T]\times\R^k\times\R^d\times\R^{dr}\longmapsto \R^d$ \ be measurable
functions.
Consider the system of semilinear  PDEs
 \vskip0.3cm\noindent  $(\mathcal{P}^{(g,F)})$ \quad $ \left\{
\begin{tabular}{ll}
 $\dfrac{\partial u(t,x)}{\partial t} + \mathcal{L}u(t,x)+F(t,x,u(t,x),
 \sigma^*\nabla
 u(t,x))=0$ & $t\in ]0,T[$, $x\in \R^k$
 \\ $u(T,x) = g(x) \quad$ $x\in\R^k$
\end{tabular} \right.$
\vskip0.1cm\noindent where $\mathcal{L}:=
\dfrac{1}{2}\displaystyle\sum_{i,j}(\sigma\sigma^*)_{ij}\partial^2_{ij}+
\displaystyle\sum_i b_i\partial_i $.

 The diffusion process associated to the operator
$\mathcal{L}$ satisfies,  $$X^{t,x}_s = x + \dint_t^s b(X^{t,x}_r) dr +  \dint_t^s
\sigma(X^{t,x}_r)dW_r, \quad\quad\; t\leq s\leq T$$

 We assume throughout this section that \ $\sigma\in\mathcal{C}^3_b(\R^k,\R^{kr})$,\ and \
$b\in\mathcal{C}^2_b(\R^k,\R^{k})$.

\vskip 0.2cm\noindent  We define,
\begin{align*}\mathcal{H}^{1+}:=
\bigcup_{\delta\geq 0,\beta>1}   \left \{v\in\mathcal{C}([0,T];
\L^\beta(\R^k,e^{-\delta\mid x\mid}dx; \R^d)): \
  \int_0^T\int_{\R^k} | \sigma^*\nabla
v(s,x)|^\beta e^{-\delta\mid x\mid} dxds <\infty \right\}
\end{align*}
\begin{definition}
 A (weak) solution of
$(\mathcal{P}^{(g,F)})$ is a function $u\in\mathcal{H}^{1+}$ such that  for
every $t\in [0,T]$ and  $ \varphi\in C_{c}^{1}([0,T]\times \R ^{d})$
\begin{align*}
\int_{t}^{T}<u(s),\frac{\partial\varphi (s)}{\partial s}> ds +<u(t),
\varphi(t)> = &  <g, \varphi (T) > + \int_t^T < F(s,.,u(s),\sigma^*\nabla
u(s)), \varphi(s)>ds
\\ &
+ \int_t^T<Lu(s),\varphi (s)>ds
\end{align*}
where $<f(s),h(s)> = \int_{\R^k}f(s,x)h(s,x)dx$. \vskip 0.2cm
\end{definition}
Observe that an integrating by part shows that,
\begin{align*}
<Lu(s),\varphi (s)> \ = \  - \int_{\R^k}\frac{1}{2}\langle \sigma^*\nabla u(s,x);
\sigma^*\nabla \varphi(s,x)\rangle dx \ - <u(s), div(\tilde{b}\varphi)(s)>
\end{align*}
 where $\tilde{b}_i := b_i
-\dfrac{1}{2}\displaystyle\sum_j\partial_j(\sigma\sigma^*)_{ij}$

%
%

\subsection{\bf{Assumptions}}
Consider the following assumptions:

 \vskip 0.1cm There
exist $\delta \geq 0 $ and $ \overline{p}>1$ such that \vskip0.4cm
\noindent\bf{(A.0)} $ g(x)\in\L^{\overline{p}}(\R^k,e^{-\delta|x|}dx;\R^d) $

\vskip0.3cm \noindent \bf{(A.1)} $F(t,x,.,.)$ is continuous $\quad$ for a.e.
$(t,x)$

\vskip0.2cm \noindent \bf{(A.2)}\; $\left\{
\begin{tabular}{ll}
 There are $\eta'\in\L^{{\frac{\overline p}{2}\vee1}}([0, T]\times
 \R^k,e^{-\delta|x|}dtdx ; \R_
 +))$, \\ $f^{0'}\in
 \L^{\overline{p}}([0, T]\times\R^k,e^{-\delta|x|}dtdx ; \R_
 +)),$ and $M, M' \in\R_+$
 such that \\ \\
$\langle y, F(t,x,y,z)\rangle \leq \eta'(t,x)+f^{0'}(t,x)|y| + (M+M'|x|)|y|^2+$
$\sqrt{M+M'|x|}|y||z|$
\end{tabular} \right.$

\vskip0.3cm \noindent \bf{(A.3)} \;  $\left\{
\begin{tabular}{ll}
 There are $ \overline{\eta}'\in\L^{q}([0, T]\times\R^k,e^{-\delta|x|}dtdx ; \R_
 +))$ (for some $q>1$), $\alpha\in ]1,\overline{p}[$ \\ and $\alpha'\in
]1,\overline{p}\wedge2[$
 such that   \\\\  $|F(t,x,y,z)|\leq \overline{\eta}'(t,x) + |y|^\alpha +
|z|^{\alpha'}$
\end{tabular} \right.$
\vskip0.4cm \noindent \bf{(A.4)}\;  $\left\{
\begin{tabular}{ll}
 There are $ K,r \in\R_+$
 such that for every
$N\in\N$  and every $ x, y,\; y',\; z,\; z' $
\\  satisfying :  \ $e^{r|x|},\;\mid
y\mid,\; \mid y'\mid,\; \mid z\mid, \;\mid z'\mid \leq N, $
\\ \\$\langle
y-y';F(t,x,y,z)-F(t,x,y',z')\rangle \leq
K\log{N}\left(\dfrac{1}{N}+|y-y'|^2\right)+$ $\sqrt{K\log{N}}|y-y'||z-z'|.$
\end{tabular} \right.$
%
%
\subsection{\bf{Existence and uniqueness for $(\mathcal{P}^{(g,F)})$ }}
\begin{theorem}\label{the51} Let $p\in]\alpha\vee\alpha',\overline{p}[$ if $M'>0$ and $p=\bar{p}$
 if $M'=0$. Under assumption \bf{(A.0)-(A.4)} we have

\vskip 0.15cm\noindent
1) The PDE $(\mathcal{P}^{(g,F)})$ has a unique (weak) solution $u$ on $[0,T]$

\vskip 0.15cm\noindent 2) For every $t\in [0,T] $ there exists
$D_t\subset\R^k$ such that \\ $i)$ $\dint_{\R^k} \1_{D_t^c}\;dx=0$ , \ where \ $D_t^c := \R^k \setminus D_t$.
\\ $ii)$ for every
$t\in [0,T]$ and every $x\in D_t$, \ the BSDE
 $(E^{(\xi^{t,x},f^{t,x})})$ has a unique solution $(Y^{t,x}, Z^{t,x})$ on $[t,T]$\\
where $\xi^{t,x}:= g(X^{t,x}_T)$ \ and \ $f^{t,x}(s,y,z):=\1_{\{s>t\}}F(s,
X^{t,x}_s, y, z)$

\vskip 0.15cm\noindent 3) For every $t\in [0,T] $
 $$ \left( u(s,X^{t,x}_s),\sigma^*\nabla
 u(s,X^{t,x}_s)\right) = \left( Y_s^{t,x},Z_s^{t,x}\right)\qquad
 \hbox{a.e.} (s,x,\omega)$$
 \vskip 0.15cm\noindent 4)  There exists a positive constant $C$  depending only on $\delta, M, M', p,
\bar{p}, $ $|\sigma|_\infty,|b|_\infty$ and $T$ such that
 $$ \begin{array}{l}\sup_{0\leq t\leq T}\dint_{\R^k}\mid u(t,x)\mid^{p}
e^{-\delta'\mid x\mid}dx+\dint_0^T\dint_{\R^k} \mid \sigma^*\nabla
 u(t,x)\mid^{p\wedge2} e^{-\delta'\mid x\mid}dtdx \\ \leq C\left(\1_{[M'\neq 0]}
+\dint_{\R^k}\mid g(x)\mid^{\overline{p}}dx+\right.  $ $\dint_{\R^k}
\dint_0^T\eta'(s,x)^{\frac{\overline{p}}{2}\vee1}ds dx+ \left. \dint_{\R^k}
\dint_0^T f^{0'}(s,x)^{\overline{p}}dsdx \right) \end{array}$$
where
$\delta'=\delta+\kappa'+\1_{[M'\neq 0]}$ and $\kappa':=
\dfrac{p\overline{p}M'T}{(\overline{p}-p)}\sup(4, \frac{2p}{p-1}).$
\end{theorem}
\vskip0.1cm

\subsection{\bf{Proof of Theorem \ref{the51}}.}
\vskip0.21cm A) \bf{Existence.}


\begin{lemma}\label{lem51}  $1)$ There exists $ \kappa
>0$ depending only on $|\sigma|_\infty,|b|_\infty$ and $T$ such that
$$\sup_{t,x} \E [\exp({\kappa \sup_{t\leq s\leq T}{\mid X^{t,x}_s-x\mid^2}})]
<\infty.\eqno{(4.0)} $$ In particular, for every $r>0$ there is a constant
$C(r,\kappa)$ such that for each $(t,x)$ $$\E [\exp({r \sup_{t\leq s\leq T}{\mid
X^{t,x}_s\mid}})]\leq C(r,\kappa) \exp({r \mid x\mid}) $$
\\
$2)$ For every $\delta \geq 0$ there exists a constant $C_{\delta,T}
>1$ such that for every $\varphi\in\L^0(\R^k),$
  $t\in [0,T]$ and $s\in [t,T]$ $$
C_{\delta,T}^{-1}\int_{\R^k}|\varphi(x)|e^{-\delta\mid x\mid}dx \leq
\E\int_{\R^k}|\varphi(X_s^{t,x})|e^{-\delta\mid x\mid}dx \leq
C_{\delta,T}\int_{\R^k}|\varphi(x)|e^{-\delta\mid x\mid}dx.\eqno{(4.2)}$$
\\
 Moreover for every $\delta \geq 0$ there
exists a constant $C_{\delta,T} >1$ such that for every
$\psi\in\L^0([0,T]\times\R^k),$
  $t\in [0,T]$ and  $s\in [t,T]$
 $$
C_{\delta,T}^{-1}\int_{\R^k}\int_t^T|\psi(s,x)|ds e^{-\delta\mid x\mid}dx \leq
\E\int_{\R^k}\int_t^T|\psi(s,X_s^{t,x})|ds e^{-\delta\mid x\mid}dx \leq
C_{\delta,T}\int_{\R^k}\int_t^T|\psi(s,x)|dse^{-\delta\mid x\mid}dx.$$
\end{lemma}


\vspace{0.5cm}\noindent {\bf{Proof.} The first assertion is well known. Its
particular case follows by using triangular and Young's inequalities. Indeed
\begin{align*}
 \E[\exp({r \sup_{t\leq s\leq T}{\mid X^{t,x}_s\mid}})]
 &\leq \exp({r \mid x\mid})\E[
\exp({r \sup_{t\leq s\leq T}{\mid X^{t,x}_s-x\mid}})]\\ &\leq \exp({r \mid
x\mid})\E[ \exp({\frac{r}{\sqrt\kappa}\sqrt\kappa \sup_{t\leq s\leq T}{\mid
X^{t,x}_s-x\mid}})]\\ &\leq \exp({\frac{r^2}{\kappa}})\exp({r \mid x\mid})\E
[\exp({\kappa \sup_{t\leq s\leq T}{\mid X^{t,x}_s-x\mid^2}})].
\end{align*}
 For the second assertion, see \cite{BM} Proposition 5.1.  \eop


\begin{lemma}\label{lem00} Let $p\in]\alpha\vee\alpha',\overline{p}[$ if $M'>0$
and $p=\bar{p}$ if $M'=0.$ Let $t\in[0,T]$.  There exists $D_t\subset\R^k$such
that \\ $i)$ $\dint_{D_t^c}1\,dx=0$ \\ $ii)$ for every $x\in D_t$
\\ $\E(\mid g(X_T^{t,x})\mid^pe^{\frac{p}{2}\int_t^T\lambda^{t,x}_sds})+
\E\left(\dint_t^T\eta'(s,X_s^{t,x})
e^{\int_t^s\lambda^{t,x}_rdr}ds\right)^\frac{p}{2}$\\
 $\indent\indent\indent\indent\indent\indent\indent\,+\E\left(\dint_t^T
f^{0'}(s,X_s^{t,x}) e^{\frac{1}{2}\int_t^s\lambda^{t,x}_rdr}ds\right)^p
+\E\dint_t^T\overline{\eta}'(s,X_s^{t,x})^{q}ds \; <+\infty$,
\\ where $ \lambda^{t,x}_s:= (M+M'|X_s^{t,x}|)\sup(4, \frac{2p}{p-1})
. $
\end{lemma}


\vskip0.1cm \noindent \bop. Using H$\ddot{o}$lder's inequality, Young's inequality and Lemma \ref{lem51} we
get
\\
$\E(\mid g(X_T^{t,x})\mid^pe^{\frac{p}{2}\int_t^T\lambda^{t,x}_sds})+
\E\left(\dint_t^T\eta'(s,X_s^{t,x})
e^{\int_t^s\lambda^{t,x}_rdr}ds\right)^\frac{p}{2}$ \\ $+\E\left(\dint_t^T
f^{0'}(s,X_s^{t,x}) e^{\frac{1}{2}\int_t^s\lambda^{t,x}_rdr}ds\right)^p
+\E\dint_t^T\overline{\eta}'(s,X_s^{t,x})^{q}ds \\ \leq $
 $ C\left(\E(\mid
g(X_T^{t,x})\mid^{\overline{p}})+ \E
\dint_t^T\eta'(s,X_s^{t,x})^{\frac{\overline{p}}{2}\vee1} + \right. $ $\left.
\E \dint_t^T f^{0'}(s,X_s^{t,x})^{\overline{p}}
+\E\dint_t^T\overline{\eta}'(s,X_s^{t,x})^{q}ds + \1_{[M'\neq 0]}e^{\kappa'
|x|}\right) $

\vskip0.1cm\noindent for some constant $C$ depending only on $M, M', p,
\bar{p}, $ $|\sigma|_\infty,|b|_\infty$ and $T$.

\vskip 0.15cm\noindent  We put, \\ $\Gamma^{t,x}:= C\big(\E(\mid
g(X_T^{t,x})\mid^{\overline{p}})+ \E
\dint_t^T\eta'(s,X_s^{t,x})^{\frac{\overline{p}}{2}\vee1} +
\E \dint_t^T f^{0'}(s,X_s^{t,x})^{\overline{p}}
+\E\dint_t^T\overline{\eta}'(s,X_s^{t,x})^{q}ds + \1_{[M'\neq 0]}e^{\kappa'
|x|}$\big).

\vskip0.2cm\noindent Using Lemma \ref{lem51}-2) and assumptions
\bf{(A.0)-(A.3)}, one can show that $$ \int_{\R^k} \Gamma^{t,x}e^{-\delta'|x|}dx
\;<\;\infty $$ where $\delta'=\delta+\kappa'+1$. The set $D_t := \{x;\quad
\Gamma^{t,x} <\infty\}$. Lemma \ref{lem00} is proved. \eop


\begin{lemma}\label{lem52}  Assume \bf{(A.0)-(A.4)}. Let
$p\in]\alpha\vee\alpha',\overline{p}[$ if $M'>0$ and $p=\bar{p}$ if $M'=0.$
Then, for every  $t\in [0,T]$ and every $x\in D_t$, \ the BSDE
 $(E^{(\xi^{t,x},f^{t,x})})$ has a unique solution $(Y^{t,x}, Z^{t,x})$
 which satisfies,  for every $t\in [0,T]$ and every
 $x\in D_t$,
$$
\begin{tabular}{ll}
  $ \E\displaystyle\big(\sup_{t\leq s\leq T}\mid Y_s^{t,x}\mid^p\big)  +
\E\big(\int_t^T \mid Z_s^{t,x}\mid^2 ds \big)^{\frac{p}{2}}$
\\ $\leq C\big[\E(\mid g(X_T^{t,x})\mid^{\overline{p}})+  $
$\E \dint_t^T\eta'(s,X_s^{t,x})^{\frac{\overline{p}}{2}\vee1}ds +  \E
\dint_t^T f^{0'}(s,X_s^{t,x})^{\overline{p}}ds + \1_{[M'\neq 0]}e^{\kappa' |x|}
\big ] $
 \end{tabular}  \eqno{(4.3)}$$ for some constant $C$  depending only on $M, M', p,
\bar{p}, $ $|\sigma|_\infty,|b|_\infty$ and $T$  .
\end{lemma}


\noindent\bf{Proof.} \ For every  $t\in[0, T]$ and $x\in D_t$, $(\xi^{t,x}, f^{t,x}) $ satisfies
\bf{(H.0)-(H.4)} with $\gamma=\inf\{\dfrac{1}{4},\dfrac{p-1}{4}\}$,
$M_s=M+M'|X^{t,x}_s|$, $K_s=\sqrt{M+M'|X^{t,x}_s|}$, $\eta_s=
\eta'(s,X_s^{t,x})$, $f_s^{0}  = f^{0'}(s,X_s^{t,x})$, $\overline{\eta}_s =
\overline{\eta}'(s,X_s^{t,x})$, $v_s = \exp({r|X^{t,x}_s|})$ and $A_N = N$.
Hence, Lemma \ref{lem52} follows from Theorem \ref{the21} and Lemma \ref{lem00}.
\eop

\vskip 0.2cm\noindent  Set,

 $g_n(x):= g(x)\1_{\{|g(x)|\leq n\}}$,
\begin{align*}
 F_{n}(t,x,y,z) \ := \ & (n^{2p}e^{|x|})^{(d+dr)}(c_1e)^2
 \1_{\{\eta'(t,x)+\overline{\eta}'(t,x)+f^{0'}(t,x)+|x|\leq n\}}
\psi(n^{-2}|y|^2)\psi(n^{-2}|z|^2) \times \\ & \ \  \int_{\R^d}\int_{\R^{dr}}
F(t,x,y-u,z-v)\Pi_{i=1}^d\psi(n^{2p}e^{|x|}u_i)\Pi_{i=1}^d
\Pi_{j=1}^r\psi(n^{2p}e^{|x|}v_{ij})dudv,
\end{align*}

$\xi^{t,x}_n := g_n(X^{t,x}_T)$ \\ and

 $f_n^{t,x}(s,y,z):=\1_{\{s>t\}}F_n(s,
X^{t,x}_s, y, z)$.

\vskip 0.2cm\noindent   It is not difficult to see that the sequence $(g_n,F_n)$ satisfies \bf{(A.0)-(A.3)}
uniformly in $n$. Hence $(\xi^{t,x}_n, f_n^{t,x})$ satisfies \bf{(H.0)-(H.3)}
uniformly in $n$. Moreover, for every $n\in \N^*$,  $(\xi^{t,x}_n, f_n^{t,x})$ is
bounded and $ f_n^{t,x}$ is globally Lipschitz.

\noindent Let $(Y^{t,x,n}, Z^{t,x,n})$ be the unique solution of BSDE
$(E^{(\xi^{t,x}_n,f_n^{t,x})})$. Let
$p\in]\alpha\vee\alpha',\overline{p}[$ if $M'>0$ and $p=\bar{p}$ if $M'=0.$ Arguing as in Lemma \ref{lem52}, we show that for every
$t,$ $x\in D_t $ and every $n\in \N^*$
$$
\begin{tabular}{ll}
  $ \E\displaystyle(\sup_{t\leq s\leq T}\mid Y_s^{t,x,n}\mid^p) +
\E\big(\dint_t^T \mid Z_s^{t,x,n}\mid^2 ds \big)^{\frac{p}{2}} \leq
C\bigg(\E\dint_t^T e^{-(\frac{\overline{p}}{2}\vee 1)|X^{t,x}_s|}ds+\E(\mid
g(X_T^{t,x})\mid^{\overline{p}})+  $
\\
\hskip 4cm $+ \E \dint_t^T\eta'(s,X_s^{t,x})^{\frac{\overline{p}}{2}\vee1}ds +
 \E \dint_t^T f^{0'}(s,X_s^{t,x})^{\overline{p}}ds + \1_{[M'\neq 0]}
e^{\kappa' |x|} \bigg) $
 \end{tabular}  \eqno{(4.4)}$$
for some constant $C=C(\bar p)$ not depending on $(t,x,n)$. To see this, use proposition
\ref{pro45} (with $h_s:=e^{-|X^{t,x}_s|}$), Proposition \ref{pro41} and the
proof of proposition \ref{pro44}-a).

\vskip 0.4cm According to \cite{BM} (see also \cite{BL}) we have
\begin{lemma}\label{lem54}
There exists a unique solution $u^n$ to the problem,

\vskip0.2cm
$(\mathcal{P}^{(g_n,F_n)})$ $ \left\{
\begin{tabular}{ll}
 $\dfrac{\partial u^n(t,x)}{\partial t} + \mathcal{L}u^n(t,x)+F_n(t,x,u^n(t,x),
  \sigma^*\nabla
 u^n(t,x))=0$, & $t\in ]0,T[$, $x\in \R^k$
 \\ $u^n(T,x) = g_n(x), \quad$ $x\in\R^k$
\end{tabular} \right.$ \vskip0.2cm \noindent
such that for every $t$ $$ u^n(s,X^{t,x}_s) = Y_s^{t,x,n} \quad \hbox{and}\quad
\sigma^*\nabla
 u^n(s,X^{t,x}_s) = Z_s^{t,x,n} \quad\quad \hbox{a.e}\;(s,\omega,x) .$$
\end{lemma}
\vskip1cm From Proposition \ref{pro44}-(ii) we have


\begin{lemma}\label{lem55}(Stability)
For every $t\in [0,T]$, $x\in D_t$ and $p'<\bar p$,

$$\lim_n \left[\E(\sup_{0\leq s\leq T}\mid Y_s^{t,x,n}-Y_s^{t,x}\mid^{p'})+
\E\left(\int_t^T \mid Z_s^{t,x,n}-Z_s^{t,x}\mid^2 ds
\right)^\frac{p'}{2}\right]=0.$$
\end{lemma}
\vskip0.5cm Using Lemma \ref{lem51}$-2)$, inequality (\ref{lem54}), Lemma \ref{lem54},
Lemma \ref{lem55} and the Lebesgue dominated convergence theorem, we obtain
\begin{lemma}\label{lem55'}(Covergence of PDE) For every $p'<\bar p$,
\begin{align*}
& \lim_{n,m}\sup_{0\leq t\leq T}\int_{\R^k}\mid u^n(t,x)-u^m(t,x)\mid^{p'}
e^{-\delta'\mid x\mid}dx =0 \\ & \lim_{n,m} \int_0^T\int_{\R^k} \mid
\sigma^*\nabla
 u^n(t,x)-\sigma^*\nabla
 u^m(t,x)\mid^{p'\wedge2} e^{-\delta'\mid x\mid} dtdx =0.
 \end{align*}
 \end{lemma}

\vskip 0.2cm
 Using Lemma \ref{lem51}, Lemma \ref{lem55'} and the fact that $\mathcal{H}^{1+}$ is
complete, we prove that  exists $u\in\mathcal{H}^{1+}$ such that for every $p'<\bar p$,

\vskip
0.2cm\noindent
  $i)$ \ $\sup_{0\leq t\leq T}\int_{\R^k}\mid u(t,x)\mid^{p'}
e^{-\delta'\mid x\mid}dx+\int_0^T\int_{\R^k} \mid \sigma^*\nabla
 u(t,x)\mid^{p'\wedge2} e^{-\delta'\mid x\mid}dtdx <\infty$

 \vskip
0.1cm\noindent
 $ii)$ \ $\lim_{n}\sup_{0\leq t\leq T}\int_{\R^k}\mid u^n(t,x)-u(t,x)\mid^{p'}
e^{-\delta'\mid x\mid}dx =0 $

 \vskip
0.1cm\noindent $iii)$ \ $ \lim_{n} \E\int_{\R^k}\left(\int_t^T
\mid \sigma^*\nabla
 u^n(s,X_s^{t,x})-\sigma^*\nabla
 u(s,X_s^{t,x})\mid^2e^{-\delta'\mid x\mid} ds\right)^\frac{p'}{2}dx =0
 \qquad \forall t\in [0,T] $

 \vskip
0.1cm\noindent
$iv)$ \quad $ \left(u(s,X^{t,x}_s), \sigma^*\nabla
 u(s,X^{t,x}_s) \right) = \left(Y_s^{t,x},Z_s^{t,x}\right) \quad \hbox{a.e.}$

 \vskip 0.5cm
In another hand, from Proposition \ref{pro42} and Proposition \ref{pro44}
we respectively have for every $t\in [0,T]$ and $x\in D_t$
$$\E\int_t^T|F_n(s,X_s^{t,x},u^n(s,X_s^{t,x}),\sigma^*\nabla
 u^n(s,X^{t,x}_s)|^{\hat{\beta}}ds\leq
 C\left(1+\Theta_{p}^{t,x,n}+
 \E\int_t^T|\overline{\eta}'(s,X_s^{t,x})|^{q}ds\right)
 $$
 and
 $$ \lim_n\E\int_t^T |F_n(s,X_s^{t,x},u^n(s,X_s^{t,x}),\sigma^*\nabla
 u^n(s,X^{t,x}_s))- F(s,X_s^{t,x},u(s,X_s^{t,x}),\sigma^*\nabla
 u(s,X^{t,x}_s)) |^{\hat{\beta}}ds =0$$

\ni where $\hat{\beta}$ is some real in $]1,\infty[$, $C$ is some  constant not depending
on $(t,x,n)$ and \\ $\Theta_{p}^{t,x,n}=\E\sup_s|Y^{t,x,n}_s|^p+
\E\big(\dint_{t}^{T}|Z^{t,x,n}_{s}|^{2} ds \big)^{\frac{p}{2}}$. \\ We deduce from Lemma \ref{lem51}, the Lebesgue dominated
convergence theorem and  inequality (4.4) that  $$
\lim_n\dint_0^T\dint_{\R^d}|F_n(s,x,u^n(s,x),\sigma^*\nabla
 u^n(s,x))- F(s,x,u(s,x),\sigma^*\nabla
 u(s,x))|^{\hat{\beta}} e^{-(1+\delta')\mid x\mid}dxds=0. $$

 \vskip0.2cm

As a consequence of Lemma \ref{lem52} and the proof of Proposition
 \ref{pro44},
we get the following existence result for the problem  $(\mathcal{P}^{(g,F)})$.


\begin{proposition}\label{pro51} Under assumptions \bf{(A.0)-(A.4)},
the PDE $(\mathcal{P}^{(g,F)})$ has a unique solution $u$ such that $
u(s,X^{t,x}_s) = Y_s^{t,x}$ and $\sigma^*\nabla
 u(s,X^{t,x}_s) = Z_s^{t,x}$. Moreover,  letting $p\in]\alpha\vee\alpha',
 \overline{p}[$ if $M'>0$ and $p=\bar{p}$ if $M'=0,$ then there is a constant
 $C$  depending only on $\delta', M, M', p, \bar{p}, $
$|\sigma|_\infty,|b|_\infty$ and $T$ such that
 $$ \begin{array}{l}\sup_{0\leq t\leq T}\dint_{\R^k}\mid u(t,x)\mid^{p}
e^{-\delta'\mid x\mid}dx+\dint_0^T\dint_{\R^k} \mid \sigma^*\nabla
 u(t,x)\mid^{p\wedge2} e^{-\delta'\mid x\mid}dtdx \\ \leq C\left(1
  +\dint_{\R^k}\mid g(x)\mid^{\overline{p}}dx+\right.

  $$\dint_{\R^k} \dint_0^T\eta'(s,x)^{\frac{\overline{p}}{2}\vee1}ds
dx+ \left. \dint_{\R^k} \dint_0^T f^{0'}(s,x)^{\overline{p}}dsdx \right)
\end{array}$$ \ni where $\delta'=\delta+\kappa'+1$ and $\kappa':=
\dfrac{p\overline{p}M'T}{(\overline{p}-p)}\sup(4, \frac{2p}{p-1}).$
\end{proposition}

\vskip0.3cm \noindent B) \bf{Uniqueness.} \vskip 0.2cm Due to the degeneracy of the diffusion coefficient, the solution of the homogeneous linear PDEs is not sufficiently smooth and hence we can not use it as a test function. In order to construct a suitable test function, we need the following lemma.   This lemma is interesting in itself since it gives a uniform estimate for a regularized degenerate PDE.

Let $\mathcal{W}_{q}^{1,2}([0, \ T]\times\R^d)$ denotes the Sobolev space of all funcions $u(t,x)$ defined on $\R_+\times\R^d$ such that both $u$ and all the generalized derivatives $D_t u$, $D_x u$, and $D^2_{xx} u$ belong to $L^q ([0, \ T]\times\R^d)$.


\begin{lemma}\label{lemmbounds}
Let $\varepsilon\in ]0,1[$, $g\in\mathcal{C}^\infty_c([0,T]\times \R^k; \R)$.
Then, the PDE

\vskip0.2cm \noindent
 $(\mathcal{P}_\varepsilon(g))$ $ \left\{
\begin{tabular}{ll}
 $\dfrac{\partial \phi^\varepsilon(t,x)}{\partial t}
 -\dfrac{1}{2}div(\sigma\sigma^*\nabla
 \phi^\varepsilon)-\varepsilon\bigtriangleup\phi^\varepsilon(t,x)
+\langle\tilde{b}(x);\nabla \phi^\varepsilon(t,x)\rangle =g(t,x) $
\\ $\phi^\varepsilon(0,x) = 0 \quad$ $x\in\R^k$
\end{tabular} \right.$
\\
has a unique solution \ $\phi^\varepsilon$ \ which satisfies :

\vskip0.15cm
$(i)$ \ \ $ \phi^\varepsilon\in\displaystyle
\bigcap_{q>\frac{3}{2}}\mathcal{W}^{1,2}_q([0,T]\times \R^k; \R)
\cap\mathcal{C}^{1,2}([0,T]\times \R^k; \R)$

\vskip0.15cm
$(ii)$ \ \ $\displaystyle \sup_{(\varepsilon,t,x)}\left\{|\dfrac{\partial \phi^\varepsilon}{\partial t
}(t,x)|+|\nabla \phi^\varepsilon (t,x) |+|\phi^\varepsilon (t,x)|\right\} <
\infty$.
\end{lemma}


\noindent\bf{Proof.} The existence and uniqueness, of the solution $\phi^\varepsilon$, follow
from \cite{LSU} (p. 318 and pp. $341-342$).
\\
We shall prove an uniform estimates
for $\phi^\varepsilon$ and for their first derivatives. These estimates  can be
established by adapting the proofs given in Krylov \cite{KR}
  pp. $330-344$. However, we give here a probabilistic
 proof which is very simple.  We assume that the
dimension $k$ is $1$. Let $X_t^\varepsilon(x)$ denotes the diffusion process
associated to the problem $(\mathcal{P}_\varepsilon(g))$. For  simplicity,
 we assume that $g$ does not depend from $t$
and the drift coefficient of $X_t^\varepsilon(x)$ is zero. The process
$X_t^\varepsilon(x)$ is then the unique (strong) solution of the following SDE

$$X_t^\varepsilon(x) = x +  \dint_0^t
\sigma_\varepsilon(X_s^\varepsilon(x))dW_s, \quad\quad\; 0\leq t\leq T$$
\\
Let $M:= \sup_{(\varepsilon,t,x)}(\vert g'(X_t^\varepsilon(x))\vert +
\vert\sigma(t,x)\vert + \vert\sigma'(t,x)\vert)$.
 Since the
coefficients $\sigma_\varepsilon$ is smooth and
 uniformly elliptic, then the solution
$\phi^\varepsilon$ belongs to $\mathcal{C}^{1,2}$. Hence, It\^o's formula shows
that,

$$\phi^\varepsilon(t,x) = -\E\int_t^T g(X_s^\varepsilon(x))ds.$$ Since
$g\in\mathcal{C}_c^\infty$, we immediately get

$$\sup_{(\varepsilon,t,x)}\left\{|\dfrac{\partial \phi^\varepsilon}{\partial t
}(t,x)|+|\phi^\varepsilon (t,x)|\right\} < \infty.$$ Since
$\sigma_\varepsilon\in\mathcal{C}^3_b$, we can show that
\begin{align*}
\vert\frac{\partial\phi^\varepsilon(t,x)}{\partial x}\vert & \leq M\E\int_t^T
\vert\frac{\partial X_s^\varepsilon(x)}{\partial x}\vert ds\\
\end{align*}
\\
It remains to show that
$\displaystyle\sup_{(\varepsilon,t,x)}\E(\vert\frac{\partial
X_t^\varepsilon(x)}{\partial x}\vert)<\infty$.
\\
Since $\vert\sigma_\varepsilon'(t,x)\vert\leq\vert\sigma'(t,x)\vert
\leq\sup_{(t,x)}\vert\sigma'(t,x)\vert \leq M$, we have
\begin{align*} \E(\vert\frac{\partial X_t^\varepsilon(x)}{\partial
x}\vert^2)&\leq 1 +
\E\int_0^t\vert\sigma_\varepsilon'(X_s^\varepsilon(x))\vert^2
\vert\frac{\partial X_s^\varepsilon(x)}{\partial x}\vert^2 ds\\ &\leq 1 +
M^2\E\int_0^t \vert\frac{\partial X_s^\varepsilon(x)}{\partial x}\vert^2 ds
\end{align*}
The Gronwall Lemma gives now the desired result. \\ In multidimensional case,
the proof can be performed similarly since it is based on the fact that the
first derivative of $\sigma_\varepsilon$ is bounded uniformly in $\varepsilon$,
which is valid in multidimensional case also, see Freidlin \cite {FR}, III $\S$ 3.2,  pp. 188-193.   Lemma
\ref{lemmbounds} is proved.
  \eop

\begin{remark}\label{estimatePhiseconde}
 \ $(i)$ \ According to Krylov estimate (because $\sigma_\varepsilon$ is uniformly elliptic), the previous proof (in dimension one) remains valid also when the coefficients $\sigma$ and $b$ are Lipschitz only.

 $(ii)$ \ Since in our situation \  $\sigma\in\mathcal{C}^3_b(\R^k,\R^{kr})$ \ and \
$b\in\mathcal{C}^2_b(\R^k,\R^{k})$, we can estimate also the second derivative of $\phi^\varepsilon$.  More precisely we have
 $$ \sup_{(\varepsilon,t,x)}\big\{\phi^\varepsilon (t,x)|+|\dfrac{\partial \phi^\varepsilon}{\partial t
}(t,x)|+|\nabla \phi^\varepsilon (t,x) |+|D^2\phi^\varepsilon (t,x)|\big\} < \infty.$$
\end{remark}

\noindent\bf{Proof of Remark \ref{estimatePhiseconde} .} \ Let $B_t$ be a $d$-dimensional Wiener process stochastically independent of $W_t$ and consider the SDE :

$$X_s^{t,x}(\varepsilon) = x + \dint_t^s \bar{b}(X_r^{t,x}(\varepsilon)) dr+  \dint_t^s
\sigma(X_r^{t,x}(\varepsilon))dW_r+\sqrt{2\varepsilon} (B_s-B_t), \quad\quad\; t\leq s\leq T$$

where $ \bar{b}(x) := \tilde{b}(x) -\dfrac{1}{2} \displaystyle\sum_j\partial_j(\sigma\sigma^*)_{.j}(x)= b(x) -
\displaystyle\sum_j\partial_j(\sigma\sigma^*)_{.j}(x)$

It\^o's formula shows that,

$$\phi^\varepsilon(T-t,x) = \E\int_t^T g(r,X_r^{t,x}(\varepsilon))dr$$

Then

$$\partial_i\phi^\varepsilon(T-t,x) = \E\int_t^T \langle\nabla g(r,X_r^{t,x}(\varepsilon));
 \partial_i X_r^{t,x}(\varepsilon)\rangle dr$$

and

$$\partial^2_{ij}
\phi^\varepsilon(T-t,x) = \E\int_t^T \langle\nabla g(r,X_r^{t,x}(\varepsilon));
 \partial^2_{ij} X_r^{t,x}(\varepsilon)\rangle +
 \langle D^2g(r,X_r^{t,x}(\varepsilon))\,\partial_{i} X_r^{t,x}(\varepsilon);
 \partial_{j} X_r^{t,x}(\varepsilon)\rangle dr$$

On other hand,
$$\partial_i(X_s^{t,x})_k(\varepsilon) = \delta_{ik} + \dint_t^s \langle\nabla\bar{b}_k(X_r^{t,x}(\varepsilon));
\partial_iX_r^{t,x}(\varepsilon)\rangle dr+ \displaystyle\sum_n \dint_t^s
\langle\nabla\sigma_{kn}(X_r^{t,x}(\varepsilon));\partial_iX_r^{t,x}(\varepsilon)\rangle dW^n_r$$
and
$$\begin{array}{ll}\partial^2_{ij}(X_s^{t,x})_k(\varepsilon) = & \dint_t^s
\langle\nabla\bar{b}_k(X_r^{t,x}(\varepsilon));
\partial^2_{ij}X_r^{t,x}(\varepsilon)\rangle dr+ \displaystyle\sum_n \dint_t^s
\langle\nabla\sigma_{kn}(X_r^{t,x}(\varepsilon));\partial^2_{ij}X_r^{t,x}(\varepsilon)\rangle dW^n_r
\\ & + \dint_t^s \langle D^2\bar{b}_k(X_r^{t,x}(\varepsilon))\partial_jX_r^{t,x}(\varepsilon);
\partial_iX_r^{t,x}(\varepsilon)\rangle dr \\ &+ \displaystyle\sum_n \dint_t^s
\langle D^2
\sigma_{kn}(X_r^{t,x}(\varepsilon))\partial_jX_r^{t,x}(\varepsilon);\partial_iX_r^{t,x}(\varepsilon)\rangle
dW^n_r\end{array}$$

It\^o's formula gives

$$\begin{array}{lll}\E|\partial_i(X_s^{t,x})_k(\varepsilon)|^4 &=& \delta_{ik}
+ 4\E\dint_t^s \langle\nabla\bar{b}_k(X_r^{t,x}(\varepsilon))\; ;\;
\partial_iX_r^{t,x}(\varepsilon)\rangle \;(\partial_i(X_r^{t,x})_k(\varepsilon))^3dr\\
 & &+ 6\displaystyle\sum_n \E\dint_t^s
|\langle\nabla\sigma_{kn}(X_r^{t,x}(\varepsilon))\; ;\;\partial_iX_r^{t,x}(\varepsilon)\rangle|^2
\;(\partial_i(X_r^{t,x})_k(\varepsilon))^2 dr
\\ & \leq & \delta_{ik} + \sup_x{(2|\nabla\bar{b}_k(x)|+ \displaystyle\sum_n |\nabla\sigma_{kn}(x)|^2)}\dint_t^s
\E|\partial_iX_r^{t,x}(\varepsilon)|^4dr

\end{array}$$

and

$$\begin{array}{lll}\E|\partial^2_{ij}(X_s^{t,x})_k(\varepsilon)|^2  & = & 2\E \dint_t^s
\langle\nabla\bar{b}_k(X_r^{t,x}(\varepsilon));
\partial^2_{ij}X_r^{t,x}(\varepsilon)\rangle \;\partial^2_{ij}(X_r^{t,x})_k(\varepsilon) dr
 \\ & &+\displaystyle\sum_n \E\dint_t^s
|\langle\nabla\sigma_{kn}(X_r^{t,x}(\varepsilon))\; ;\;\partial^2_{ij}X_r^{t,x}(\varepsilon)\rangle|^2 dr
\\ & &+ 2\E\dint_t^s \langle D^2\bar{b}_k(X_r^{t,x}(\varepsilon))\;\partial_jX_r^{t,x}(\varepsilon)\; ;\;
\partial_iX_r^{t,x}(\varepsilon)\rangle \;\partial^2_{ij}(X_r^{t,x})_k(\varepsilon)dr \\
& &+ \displaystyle\sum_n \E\dint_t^s |\langle D^2
\sigma_{kn}(X_r^{t,x}(\varepsilon))\;\partial_jX_r^{t,x}(\varepsilon)\;
;\;\partial_iX_r^{t,x}(\varepsilon)\rangle|^2 dr\\ & \leq  & \sup_x{(2|\nabla\bar{b}_k(x)|+2|D^2\bar{b}_k(x)|+
\displaystyle\sum_n |\nabla\sigma_{kn}(x)|^2)}\E \dint_t^s |
\partial^2_{ij}X_r^{t,x}(\varepsilon)|^2 dr
 \\ & &+ \sup_x{(|D^2\bar{b}_k(x)|+\displaystyle\sum_n |D^2
\sigma_{kn}(x)|^2)}\dint_t^s  \E|\partial_jX_r^{t,x}(\varepsilon)|^4+ \E|\partial_iX_r^{t,x}(\varepsilon)|^4dr
\end{array}$$

We deduce that

$$\begin{array}{lll}\E|\partial_i(X_s^{t,x})(\varepsilon)|^4  & \leq & k^2 +
k^2\sum_j\sup_x{(2|\nabla\bar{b}_n(x)|+ \displaystyle\sum_n |\nabla\sigma_{jn}(x)|^2)}\dint_t^s
\E|\partial_iX_r^{t,x}(\varepsilon)|^4dr\\ & \leq  & k^2 e^{k^2T\sum_j\sup_x{(2|\nabla\bar{b}_n(x)|+ \sum_n
|\nabla\sigma_{jn}(x)|^2)}}\qquad \hbox{(Gronwall's Lemma )}

\end{array}$$

and

$$\begin{array}{lll}\E|\partial^2_{ij}(X_s^{t,x})(\varepsilon)|^2
 & \leq  & k\sup_x{(2|\nabla\bar{b}_k(x)|+2|D^2\bar{b}_k(x)|+
\displaystyle\sum_n |\nabla\sigma_{kn}(x)|^2)}\E \dint_t^s |
\partial^2_{ij}X_r^{t,x}(\varepsilon)|^2 dr
 \\ & &+ k^3T\sup_x{(|D^2\bar{b}_k(x)|+\displaystyle\sum_n |D^2
\sigma_{kn}(x)|^2)}k^2 e^{k^2T\sum_j\sup_x{(2|\nabla\bar{b}_n(x)|+ \sum_n |\nabla\sigma_{jn}(x)|^2)}}\\
& \leq & k^3T\sup_x{(|D^2\bar{b}_k(x)|+\displaystyle\sum_n |D^2 \sigma_{kn}(x)|^2)}k^2
e^{k^2T\sum_j\sup_x{(2|\nabla\bar{b}_n(x)|+ \sum_n |\nabla\sigma_{jn}(x)|^2)}}\\ & & \times\;
e^{kT\sup_x{(2|\nabla\bar{b}_k(x)|+2|D^2\bar{b}_k(x)|+ \sum_n |\nabla\sigma_{kn}(x)|^2)}} \qquad \hbox{(Gronwall's
Lemma )}
\end{array}$$

Since $g\in\mathcal{C}_c^\infty$, $\sigma\in\mathcal{C}^3_b(\R^k,\R^{kr})$ and $b\in\mathcal{C}^2_b(\R^k,\R^{k})$
we  get

$$ \sup_{(\varepsilon,t,x)}\left\{\phi^\varepsilon (t,x)|+|\dfrac{\partial \phi^\varepsilon}{\partial t
}(t,x)|+|\nabla \phi^\varepsilon (t,x) |+|D^2\phi^\varepsilon (t,x)|\right\} < \infty.$$ Lemma \ref{lemmbounds} is
proved.
  \eop


\begin{lemma}\label{lemm1}
$0$ is the unique solution of the PDE \vskip0.2cm \noindent
 $(\mathcal{P}^{(0,-div(\tilde{b})(x)y)})$ $ \left\{
\begin{tabular}{ll}
 $\dfrac{\partial w(t,x)}{\partial t} + \mathcal{L}w(t,x) +
 div(\tilde{b})(x)w(t,x)=0$
 & $t\in ]0,T[$, $x\in \R^k$
 \\ $w(T,x) = 0 \quad$ $x\in\R^k$
\end{tabular} \right.$

\vskip0.2cm \noindent satisfying for some $\beta >1$
\begin{equation}\label{equa1}
\sup_{0\leq t\leq T}\int_{\R^k}\mid w(t,x)\mid^{\beta}+\mid w(t,x)\mid
dx+\int_{0}^T\int_{\R^k} \mid \sigma^*\nabla
 w(t,x)\mid^{\beta}+\mid \sigma^*\nabla
 w(t,x)\mid dtdx <\infty.
\end{equation}
\end{lemma}


 \noindent\bf{Proof.} Let $w$ be a solution of
$(\mathcal{P}{(0,-div(\tilde{b})(x)y)})$ satisfying (\ref{equa1}) and consider
$w_n\in\mathcal{C}^\infty_c(R^k)$ such that
$$\displaystyle\int_{0}^T\displaystyle\int_{\R^k} |
  w(s,x)-w_n(s,x)| dxds +
\displaystyle\int_{0}^T\displaystyle\int_{\R^k} | \sigma^*\nabla
(w(s,x)-w_n(s,x))| dxds \rightarrow 0.  $$

\vskip 0.3cm\noindent Let $\varepsilon\in ]0,1[$,
$g\in\mathcal{C}^\infty_c([0,T]\times \R^k; \R)$ and consider the unique
solution  $\phi^\varepsilon\in\displaystyle
\cap_{q>\frac{3}{2}}\mathcal{W}^{1,2}_q([0,T]\times \R^k; \R)
\cap\mathcal{C}^{1,2}([0,T]\times \R^k; \R)$ of the following problem

\vskip0.2cm \noindent
 $(\mathcal{P}_\varepsilon(g))$ $ \left\{
\begin{tabular}{ll}
 $\dfrac{\partial \phi^\varepsilon(t,x)}{\partial t}
 -\dfrac{1}{2}div(\sigma\sigma^*\nabla
 \phi^\varepsilon)-\varepsilon\bigtriangleup\phi^\varepsilon(t,x)
+\langle\tilde{b}(x);\nabla \phi^\varepsilon(t,x)\rangle =g(t,x)\rangle $
\\ $\phi^\varepsilon(0,x) = 0 \quad$ $x\in\R^k$
\end{tabular} \right.$

\vskip0.2cm \noindent  The existence and uniqueness of $\phi^\varepsilon$
follows from Lemma \ref{lemmbounds}.

\vskip 0.3cm\noindent
 Let
$(\psi_i)_{i\in\N}\subset\mathcal{C}^\infty_c(R^k)$ be such that $\psi_i\in
[0,1]$,\quad $\psi_i\rightarrow 1$ uniformly on every compact set and \quad
$\nabla\psi_i\rightarrow 0$ uniformly on $\R^k. $ By considering $\phi^\varepsilon\psi_i$ as a test function, we have

\begin{align*}\int_{0}^T\int_{\R^k} & \left[w\dfrac{\partial
\phi^\varepsilon}{\partial t }+\dfrac{1}{2}\langle\sigma^*\nabla
w;\sigma^*\nabla \phi^\varepsilon \rangle+w\langle\tilde{b};\nabla
\phi^\varepsilon \rangle \right]\psi_i dxdt +\\ &
\int_{0}^T\int_{\R^k}\left[\dfrac{1}{2}\langle\sigma^*\nabla w;\sigma^*\nabla
\psi_i\rangle + w\langle\tilde{b};\nabla
\psi_i\rangle\right]\phi^\varepsilon dxdt =0.
\end{align*}
Introducing $w_n$ and integrating by part we obtain
$$ \int_{0}^T\int_{\R^k}
w_n\psi_i\left[\dfrac{\partial \phi^\varepsilon}{\partial t
}-\dfrac{1}{2}div(\sigma\sigma^*\nabla \phi^\varepsilon)
+\langle\tilde{b};\nabla \phi^\varepsilon \rangle \right]dtdx =
\chi_1^{\varepsilon,i}(n)+\chi_2^{\varepsilon,n}(i),
$$
where
\begin{align*} \chi_1^{\varepsilon,i}(n):= -\int_{0}^T  \int_{\R^k}
\left[(w-w_n)\dfrac{\partial \phi^\varepsilon}{\partial t
}+\dfrac{1}{2}\langle\sigma^*\nabla (w-w_n);\sigma^*\nabla \phi^\varepsilon
\rangle+ \right. \left. (w-w_n)\langle\tilde{b};\nabla \phi^\varepsilon
\rangle \right]\psi_i dxdt
\end{align*}
and $$\chi_2^{\varepsilon,n}(i):=
-\int_{0}^T\int_{\R^k}\langle\dfrac{1}{2}\phi^\varepsilon \sigma\sigma^*\nabla
w+\phi^\varepsilon w\tilde{b}-\dfrac{1}{2}w_n
\sigma\sigma^*\nabla\phi^\varepsilon \; ;\;\nabla \psi_i\rangle dxdt. $$ From
Lemma \ref{lemmbounds}, we have $$ \sup_\varepsilon
\sup_{(t,x)}\left\{|\dfrac{\partial \phi^\varepsilon}{\partial t
}(t,x)|+|\nabla \phi^\varepsilon (t,x) |+|\phi^\varepsilon (t,x)|\right\} <
\infty.$$
\\
Hence  $$
\sup_{\varepsilon,i}|\chi_1^{\varepsilon,i}(n)|\longrightarrow_{n\rightarrow\infty}
0 $$ and $$\sup_{\varepsilon,n}|\chi_2^{\varepsilon,n}(i)|
\longrightarrow_{i\rightarrow\infty} 0.$$ Observe that an integrating by part shows that \
$\int_{0}^T\int_{\R^k} w_n\psi_i \bigtriangleup\phi^\varepsilon dxdt = -
\int_{0}^T\int_{\R^k} \nabla (w_n\psi_i) \nabla\phi^\varepsilon dxdt$, \ then
use the Lebesgue dominated convergence theorem to deduce that
\begin{align*}\int_{0}^T\int_{\R^k} w g(t,x)dxdt & =
\lim_n\lim_i\lim_\varepsilon \int_{0}^T\int_{\R^k} w_n\psi_i
(g(t,x)+\varepsilon\bigtriangleup\phi^\varepsilon)dxdt
\\ &
=\lim_n\lim_i\lim_\varepsilon
(\chi_1^{\varepsilon,i}(n)+\chi_2^{\varepsilon,n}(i)) \\
& =0.
\end{align*}
Lemma
\ref{lemm1} is proved. \eop

\vskip 0.15cm
\bop\bf{of uniqueness for $(\mathcal{P}^{(g,F)})$.} The proof is divided into
three steps.
\vskip 0.15cm \noindent\bf{Step1.} \it {$0$ is the unique solution of $(\mathcal{P}^{(0,0)})$
satisfying the inequality (\ref{equa1})} Lemma 4.8.

 \vskip 0.15cm Let $w_1$ be a solution of
$(\mathcal{P}^{(0,0)})$  satisfying the inequality (\ref{equa1}) Lemma 4.8.
 Then, by Lemma \ref{lemm1} it is also the unique solution of
 $(\mathcal{P}^{(0,div\tilde{b}(x)y-div\tilde{b}(x)w_1(t,x))})$ satisfying the inequality
 (\ref{equa1}) Lemma 4.8. Indeed, if $u$ is a solution of
 $(\mathcal{P}^{(0,div\tilde{b}(x)y-div\tilde{b}(x)w_1(t,x))})$,
 then $u-w_1$ is a solution of $(\mathcal{P}^{(0,div\tilde{b}(x)y)})$ and
 hence $u-w_1 = 0$ by Lemma \ref{lemm1}. \\
 From Proposition \ref{pro51}, the process
  $(w_1(s,X^{t,x}_s),\sigma^*\nabla w_1(s,X^{t,x}_s))$
 is the unique solution of
  BSDE \\  $(E^{(0,div\tilde{b}(X^{t,x}_s)y-div\tilde{b}(X^{t,x}_s)w_1(s,X^{t,x}_s)})$.
  Thanks to the uniqueness of this BSDE and Lemma \ref{lem51}-2), we get $w_1=0$.
  \\\\\bf{Step2.} \it{ $0$ is the unique solution of
  $(\mathcal{P}^{(0,0)})$.}

   \vskip 0.15cm Let $w_1$ be a  solution of
 $(\mathcal{P}^{(0,0)})$. Since $w_1\in\mathcal{H}^{1+}$, then there
 exist
  $\beta'>1,\delta'\geq 0$ such that,

 $$ \sup_{0\leq t\leq T}\int_{\R^k}\mid w_1(t,x)\mid^{\beta'}
e^{-\delta'\mid x\mid}dx+\int_{0}^T\int_{\R^k} \mid \sigma^*\nabla
 w_1(t,x)\mid^{\beta'} e^{-\delta'\mid x\mid}dxdt <\infty.$$

 \noindent Let $\delta > \delta'$ and set $\tilde{w}_1 := w_1f(x)$
 where $f\in \mathcal{C}^2(\R^k;\R^*_+)$ such that
  $f(x)=e^{-\delta |x|}$ if $\mid x\mid > 1$.

 \noindent By Lemma \ref{lemm1}, $\tilde{w}_1$ is the unique
 solution to the PDE
 \vskip0.5cm\noindent
 $(\mathcal{P}_1^{(0,0)})$
 \hskip 1cm $ \left\{
\begin{tabular}{ll}
 $\dfrac{\partial w(t,x)}{\partial t} + \mathcal{L}w(t,x) +
 div(\tilde{b})(x)w(t,x)
 + H(x)
 \tilde{w}_1(t,x)+
  \langle\overline{H}(x),\sigma^*\nabla\tilde{w}_1(t,x)\rangle=0$
  \\ $w(T,x) = 0 $
\end{tabular} \right.$ \vskip0.2cm \noindent satisfying the inequality
(\ref{equa1}) Lemma 4.8, where $H$ and $\overline{H}$ are some
  bounded and continuous functions.

\noindent Proposition \ref{pro51} implies that
$(\tilde{w}_1(s,X^{t,x}_s),\sigma^*\nabla\tilde{w}_1(s,X^{t,x}_s))$
 is the unique solution of the BSDE \\ $(E^{(0,div(\tilde{b})(X^{t,x}_s)y + H(X^{t,x}_s)
 \tilde{w}_1(s,X^{t,x}_s)+
\langle\overline{H}(X^{t,x}_s),\sigma^*\nabla\tilde{w}_1(s,X^{t,x}_s)
\rangle)})$.
   Hence $\tilde{w}_1=0 $,  which implies that  $w_1=0$.

 \vskip 0.2cm\noindent \bf{Step 3.}  \it{$(\mathcal{P}^{(g,F)})$
 has a unique
  solution if and only if $0$ is the unique solution of
  $(\mathcal{P}^{(0,0)})$}.

  \vskip 0.15cm By Proposition \ref{pro51}, there exists a unique solution $u$
  of the problem $(\mathcal{P}^{(g,F)})$ such that,  $u(s,X^{t,x}_s) =
  Y_s^{t,x}$ \
  and \ $\sigma^*\nabla
 u(s,X^{t,x}_s) = Z_s^{t,x}$.
 \\ Let $u'$ be another solution of
  $(\mathcal{P}^{(g,F)})$
   and set
$$\hat{F}(t,x)=F(s,x,u(s,x),\sigma^*\nabla
 u(s,x))-F(s,x,u'(s,x),\sigma^*\nabla
 u'(s,x)). $$
 The function \ $w:=u-u'$ is then a solution of the problem
 \vskip0.2cm \noindent
 $(\mathcal{P}^{(0,\hat{F})})$ \qquad $ \left\{
\begin{tabular}{ll}
 $\dfrac{\partial w(t,x)}{\partial t} + \mathcal{L}w(t,x)+\hat{F}
 (t,x)=0$
 & $t\in ]0,T[$, $x\in \R^k$
 \\ $w(T,x) = 0 \quad$ $x\in\R^k$
\end{tabular} \right.$

\vskip0.2cm\noindent
In  other hand, since $(0,\hat{F})$ satisfies assumptions \bf{(A.0)-(A.4)}, then
Proposition \ref{pro51} ensures the existence of a unique solution $\hat{w}$
to the problem
 $(\mathcal{P}^{(0,\hat{F})})$ such that, $
\hat{w}(s,X^{t,x}_s) = \hat{Y}_s^{t,x}$ \ and \ $\sigma^*\nabla
 \hat{w}(s,X^{t,x}_s) = \hat{Z}_s^{t,x}$, where \
 ($\hat{Y}_s^{t,x}, \hat{Z}_s^{t,x})$ is the unique solution
 of $$ \hat{Y}_s^{t,x} = \int_s^T \hat{F}(r,X^{t,x}_r)dr -
 \int_s^T\hat{Z}_r^{t,x}dW_r $$
The uniqueness of $(\mathcal{P}^{(0,\hat{F})})$ (which follows from step 2)
allows us to deduce that $$ u'(s,X^{t,x}_s)= Y_s^{t,x}-\hat{Y}_s^{t,x}
\quad\quad\hbox{and}\quad\quad\sigma^*\nabla
 u'(s,X^{t,x}_s) = Z_s^{t,x}-\hat{Z}_s^{t,x}. $$
 This implies that $u'(t,X_s^{t,x})$ is a solution to  BSDE $(E^{(g,F)})$.
 The uniqueness of this BSDE shows that $u'(t,X_s^{t,x}) =
 u(t,X_s^{t,x})$. We get that $u(t,x) = u'(t,x)$ $a.e.$ by using
  Lemma \ref{lem51}-2).
Theorem \ref{the51} is proved. \eop


\vskip 0.3cm \noindent As consequence, we have :
 \noindent Let $ g \in \L^{p}([0,
T]\times\R^k,e^{-\delta|x|}dx ; \R^d)$ for some $p>1$ and $\delta\geq 0.$ Let
$A : [0,T]\times \R^k \longrightarrow \R^{d\times d }$, $B  : [0,T]\times \R^k
\longrightarrow (\R^d)^{dr}$ and $C : [0,T]\times \R^k \longrightarrow
\R^{d\times d }$ be measurable functions which satisfy :
\par
 There exists a positive constant $K$ such that for
all $(t,x)$ $$ \|A(t,x)\| + \|B(t,x)\|^2 \leq K (1+|x|),\;\;\|C(t,x)\|\leq K
\;\;\hbox{and}\;\; C(t,x) \geq 0.$$
\par\noindent We  then have
\begin{proposition}\label{der}
Let $ g \in \L^{p}([0,
T]\times\R^k,e^{-\delta|x|}dx ; \R^d)$ for some $p>1$ and $\delta\geq 0.$ Let
$A : [0,T]\times \R^k \longmapsto \R^{d\times d }$, $B  : [0,T]\times \R^k
\longmapsto (\R^d)^{dr}$ and $C : [0,T]\times \R^k  \longmapsto
\R^{d\times d }$ be measurable functions. Assume that there exists a positive constant $K>0$ such that for
every \  $(t,x)$, \ $0\leq C(t,x) \leq K$ and,
\par
  $ \|A(t,x)\| + \|B(t,x)\|^2 \leq K (1+|x|),\;\;\|C(t,x)\|\leq K $
\par
Then, the PDE \vskip0.2cm \noindent
  $ \left\{
\begin{tabular}{l}
 $\dfrac{\partial w(t,x)}{\partial t} + \mathcal{L}w(t,x) +
 A(t,x)w(t,x)+\langle\langle\; B(t,x);\,\sigma^*\nabla w(t,x)\;\rangle\rangle
-C(t,x)w(t,x)\log|w(t,x)|=0,$
 \\ $w(T,x) = g(x) \quad$ $x\in\R^k$
\end{tabular} \right.$
\vskip0.2cm\noindent has a unique solution $w$
and $(w(s,X^{t,x}_s),\sigma^*\nabla w(s,X^{t,x}_s))$ is the unique solution of
$$ E^{\big( g(X^{t,x}_T),\;A(s,X^{t,x}_s)y \,+ \langle\langle
B(s,X^{t,x}_s);z\rangle\rangle\, - C(s,X^{t,x}_s)y\log|y|\big)},$$ \ni where $
\langle\langle B ;z\rangle\rangle :=
\displaystyle\sum_{i=1}^d\displaystyle\sum_{j=1}^r B_{ij}Z_{ij}.$
\end{proposition}

\noindent  Set \ \ $ F(t,x,y,z):= A(t,x)y \,+ \langle\langle
B(t,x);z\rangle\rangle\, - C(t,x)y\log|y|. $ \\
Arguing as in the
introductory examples, we show the following claims 1)--3). The
claim 2) follows by using Young's inequality.
 \\
\ni $1)\; \langle y, F(t,x,y,z)\rangle \leq K +
(K+K|x|)|y|^2+\sqrt{K+K|x|}|y||z|$
\\ \\
\ni $2)$ for all  $ \varepsilon>0$ there is a constant $C_\varepsilon$ such
that \\
\\
$|F(t,x,y,z)|\leq C_\varepsilon(1+|x|^{C_\varepsilon} + |y|^{1+\varepsilon} +
|z|^{1+\varepsilon})$
\\ \\
\ni $3)$  for every $N>3$  and every $ x, y,\; y'\; z,\; z' $ satisfying \ \ $
e^{|x|},\;\mid y\mid,\; \mid y'\mid,\; \mid z\mid, \;\mid z'\mid \leq N : $
\\ \\ $\langle
y-y';F(t,x,y,z)-F(t,x,y',z')\rangle \leq
K'\log{N}\left(\dfrac{1}{N}+|y-y'|^2\right)+$ $\sqrt{K'\log{N}}|y-y'||z-z'|,$

\ni where $ K':=1+4Kd+K^2.$
\vskip 0.15cm \ni So assumptions \bf{(A.0)-(A.4)} are satisfied for $(g,F).$ \eop
\vskip 1cm
\noindent\bf{Acknowledgments.} We are grateful to Etienne Pardoux for various discussions about this work.

\vskip 0.3cm\noindent\bf{Supports.}
  \ The first author is partially supported by CMEP, PHC Tassili n$^{\circ}$ 07MDU705 and Marie Curie ITN n°. 213841-2.
\\ The second author is partially supported by
Ministerio de Educacion y Ciencia, grant number SB2003-0117 and
CMIFM, PHC Volubilis n$^{\circ}$ MA/06/142. He  would like to thank
the "Centre de Recerca Matem\`atica" for their kind hospitality.
\\ The third author is partially supported by CMIFM, PHC Volubilis n$^{\circ}$ MA/06/142 and Marie Curie ITN n°. 213841-2.


%
%

\end{document}